\documentclass{amsart}

\usepackage{graphics}
\usepackage[usenames]{color}
\usepackage{graphicx}
\usepackage{epsfig}
\usepackage{float}
\usepackage{subfigure}
\usepackage{enumerate}
\usepackage{amssymb,amsmath}
 \usepackage{verbatim}
\usepackage{lineno}
\usepackage[all]{xy}

\newtheorem{my_thm}{Theorem}[section]
\newtheorem{my_lem}[my_thm]{Lemma}
\newtheorem{my_prop}[my_thm]{Proposition}
\newcommand{\norm}[1]{\left\Arrowvert \, #1 \, \right\Arrowvert}

\numberwithin{equation}{section}

\numberwithin{figure}{section}

\numberwithin{table}{section}

\begin{document}
\title{Bounded Domain Problem for the Modified Buckley-Leverett Equation}

\author{Ying Wang}
\address{Department of Mathematics, The Ohio State University, 231 West 18th Ave, Columbus, 
OH 43210}
\curraddr{School of Mathematics, University of Minnesota, 127 Vincent Hall 206 Church St SE, Minneapolis, MN 55455}
\email{wang@math.umn.edu}

\author{Chiu-Yen Kao}
\address{Department of Mathematics, The Ohio State University, 231 West 18th Ave, Columbus, 
OH 43210; Mathematical Biosciences Institute, The Ohio State University, 1735 Neil Ave, Columbus, OH 43210}
\email{kao@math.ohio-state.edu}
\thanks{This work is supported in part by NSF Grant DMS-0811003 and an Alfred P. Sloan Fellowship.}

 
\linenumbers
\renewcommand{\thefootnote}{\fnsymbol{footnote}}


\renewcommand{\thefootnote}{\arabic{footnote}}


\date\today


\maketitle

\begin{abstract}
The focus of the present study is the modified  Buckley-Leverett
(MBL) equation describing two-phase flow in porous media. The MBL
equation differs from the classical Buckley-Leverett (BL) equation
by including a balanced diffusive-dispersive combination. The
dispersive term is a third order mixed derivatives term, which
models the dynamic effects in the pressure difference between the
two phases. The classical BL equation gives a monotone water
saturation profile for any Riemann problem;  on the contrast, when
the dispersive parameter is large enough, the MBL equation
delivers non-monotone water saturation profile for certain Riemann
problems as suggested by the experimental observations. In this
paper, we first show that the solution of the finite interval
$[0,L]$ boundary value problem converges to that of the half-line
$[0,+\infty)$ boundary value problem for the MBL equation as
$L\rightarrow +\infty$. This result provides a justification for
the use of the finite interval boundary value problem in numerical
studies for the half line problem. Furthermore, we extend the
classical central schemes for the hyperbolic conservation laws to
solve the MBL equation which is of pseudo-parabolic type.
Numerical results confirm the existence  of non-monotone water
saturation profiles consisting of constant states separated by
shocks.
\end{abstract}

%
%


\pagestyle{myheadings}
\thispagestyle{plain}

\section{Introduction}
\label{intro}
\begin{figure}[htbp]
\subfiguretopcaptrue
\begin{center}
\subfigure[$f(u)=\frac{u^2}{u^2+M(1-u)^2}$]{\label{f_figure}\includegraphics[width=0.4\textwidth]{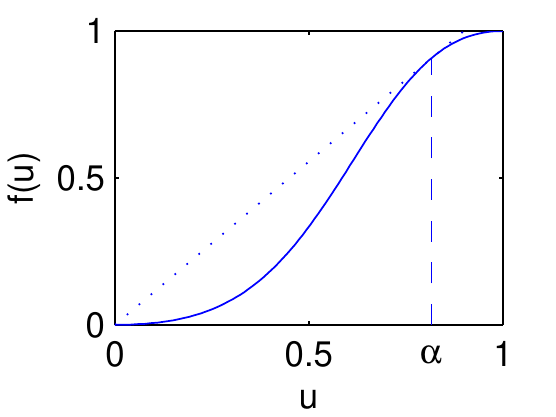}}
\subfigure[$f'(u)=\frac{2Mu(1-u)}{(u^2+M(1-u)^2)^2}$]{\label{fprime_figure}\includegraphics[width=0.4\textwidth]{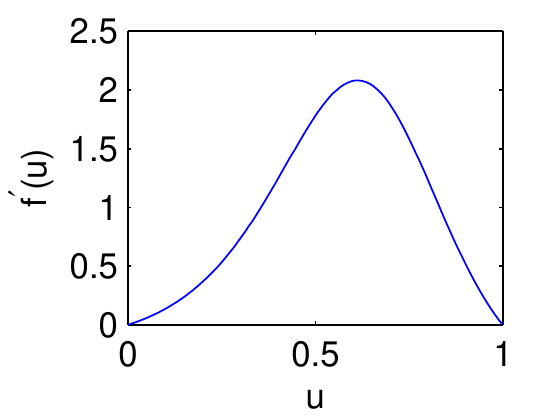}}
\vspace*{-0.5cm}\caption{
$f(u)$ and $f'(u)$ with $M=2$.} 
\label{fandfprime}
\end{center}
\end{figure}

The classical Buckley-Leverett (BL) equation \cite{Buckley} is a simple model for two-phase fluid flow in a porous medium. One application is secondary recovery by water-drive in oil reservoir simulation. In one space dimension the equation has the standard conservation form
\begin{eqnarray}
u_t+(f(u))_x =0  \qquad&\mathrm{in} & \qquad Q=\{(x,t): x>0, t>0\}\nonumber\\
u(x,0) = 0 \qquad &&\qquad x\in(0,\infty)\label{BL}\\
u(0,t) =u_B \qquad &&\qquad  t\in[0,\infty) \nonumber
\end{eqnarray}
with the flux function $f(u)$ being defined as
\begin{eqnarray}
f(u) = \left\{
\begin{array}{ll}
0 & u<0 ,\\
 \frac{u^2}{u^2+M(1-u)^2} & 0\le u \le 1 ,\\
 1 & u>1 .
\end{array}
 \right.\label{flux}
\end{eqnarray}
In this content, $u:\bar{Q}\rightarrow[0,1]$ denotes the water saturation (e.g. $u=1$ means pure water, and $u=0$ means pure oil), $u_B$ is a constant which indicates water saturation at $x=0$, and $M>0$ is the water/oil viscosity ratio. The classical BL equation (\ref{BL}) is a prototype for conservation laws with convex-concave flux functions. 
The graph of $f(u)$ and $f'(u)$ with $M=2$ is given in Figure \ref{fandfprime}.

Due to the possibility of the existence of shocks in the solution of the hyperbolic conservation laws (\ref{BL}), the weak solutions are sought. The function $u\in L^\infty(Q)$ is called a weak solution of the conservation laws (\ref{BL}) if
\begin{equation*}
\int_Q\left\{u\frac{\partial\phi}{\partial
t}+f(u)\frac{\partial\phi}{\partial x}\right\} = 0\qquad
\text{for all }\quad \phi\in C_0^\infty(Q) .
\end{equation*}
Notice that the weak solution is not unique. Among the weak
solutions, the entropy solution is physically relevant and unique.
The weak solution that satisfies Oleinik entropy condition \cite{Oleinik}
\begin{eqnarray}
\frac{f(u)-f(u_l)}{u-u_l} \ge  s  \ge \frac{f(u)-f(u_r)}{u-u_r} \qquad \mathrm{for~ all~ } u 
\mathrm{~ between ~} u_l \mathrm{ ~and ~} u_r
\label{oleinik}
\end{eqnarray}
is the entropy solution, where $u_l$, $u_r$ are the function values to the left and right of the shock respectively, and the shock speed $s$  satisfies Rankine-Hugoniot jump condition  \cite{Rankine, Hugoniot}
\begin{equation}
s=\frac{f(u_l)-f(u_r)}{u_l-u_r} .
\label{s}
\end{equation}

\begin{figure}[htbp]
\subfiguretopcaptrue
\begin{center}
\subfigure[$u_B = 0.7$]{\label{BL1}\includegraphics[width=0.4\textwidth]{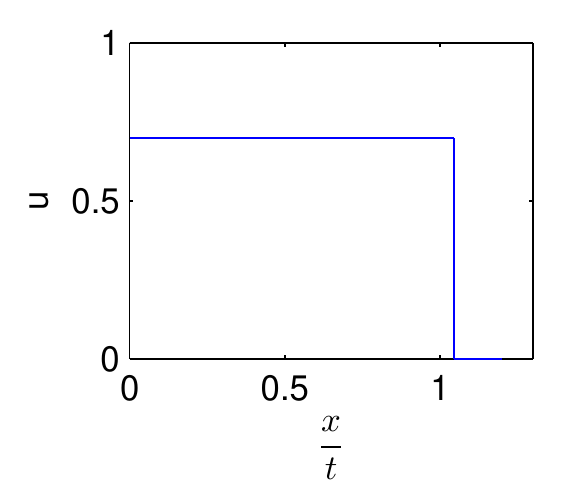}}
\subfigure[$u_B = 0.98$]{\label{BL2}\includegraphics[width=0.4\textwidth]{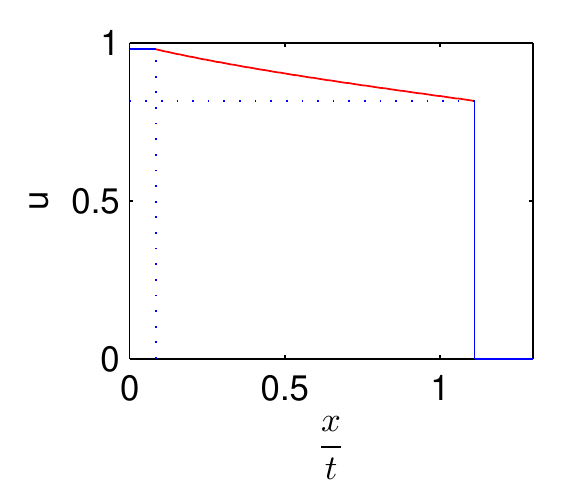}}
\vspace{-0.5cm}\caption{The entropy solution of the classical BL equation ($M=2,~
\alpha=\sqrt{\frac{2}{3}}\approx 0.8165$).  \subref{BL1}
$0<uB=0.7\le\alpha$, the solution consists of one shock at
$\frac{x}{t}=\frac{f(u_B)}{u_B}$;   \subref{BL2}
$\alpha<uB=0.98<1$, the solution consists of  a rarefaction
between $u_B$ and $\alpha$ for $f'(u_B)<\frac{x}{t}<f'(\alpha)$
and  a shock at $\frac{x}{t}=\frac{f(\alpha)}{\alpha}$.}
\label{BLfig}
\end{center}
\end{figure}

The classical BL equation (\ref{BL}) with flux function $f(u)$ as
given in (\ref{flux}) has been well studied (see \cite{Leveque}
for an introduction). Let $\alpha$ be the solution of
$f'(u)=\frac{f(u)}{u}$, i.e.,
\begin{equation}
\alpha=\sqrt{\frac{M}{M+1}}.
\label{alpha}
\end{equation}
The entropy solution of the classical BL equation can be classified into two categories:
\begin{enumerate}
 \item If $0<u_B\le\alpha$, the entropy solution has a single shock at $\frac{x}{t}=\frac{f(u_B)}{u_B}$.
 \item If $\alpha<u_B<1$, the entropy solution contains a rarefaction between $u_B$ and $\alpha$ for $f'(u_B)<\frac{x}{t}<f'(\alpha)$ and  a shock at $\frac{x}{t}=\frac{f(\alpha)}{\alpha}$.
\end{enumerate}
These two types of solutions are shown in Figure \ref{BLfig} for $M=2$.
In either case, the entropy solution  of the classical BL equation (\ref{BL}) is a non-increasing function of $x$ at any given time $t>0$.
However, the experiments of two-phase flow in porous medium reveal complex infiltration profiles, which may involve overshoot, i.e., profiles may not be monotone \cite{Dicarlo}. This suggests the need of modification to the classical BL equation (\ref{BL}).

To better describe the infiltration profiles, we go back to the origins of (\ref{BL}). Let $S_i$ be the saturation of water/oil ($i=w,o$) and assume that the medium is completely saturated, i.e. $S_w+S_o=1$. The conservation of mass gives
\begin{equation}
\phi\frac{\partial S_i}{\partial t}+\frac{\partial q_i}{\partial x} = 0 \label{governing0}
\end{equation}
where $\phi$ is the porosity of the medium (relative volume occupied by the pores) and $q_i$ denotes the discharge of 
water/oil with $q_w+q_o=q$, which is assumed to be a constant in space due to the complete saturation assumption. 
Throughout of this work, we consider it constant in time as well.
By Darcy's law 
\begin{equation}
q_i=-k\frac{{k_r}_i(S_i)}{\mu_i}\frac{\partial P_i}{\partial x}\label{darcy},\qquad i=w,o
\end{equation}
where $k$ denotes the absolute permeability, ${k_r}_i $ is the relative permeability and $\mu_i$ is the viscosity of water/oil.
Instead of considering constant capillary pressure as adopted by the classical BL equation (\ref{BL}),
Hassanizadeh and Gray \cite{Hassanizadeh, Hassanizadeh1} have defined the dynamic capillary pressure as
\begin{equation}
P_c =  P_o-P_w = p_c(S_w)-\phi\tau\frac{\partial S_w}{\partial t}
\label{capillary}
\end{equation}
where $p_c(S_w)$ is the \emph{static} capillary pressure and
$\tau$ is a positive constant, and $\frac{\partial S_w}{\partial
t}$ is the dynamic effects. Using Corey \cite{Corey, duijn2}
expressions with exponent 2, ${k_r}_w(S_w) = S_w^2, ~{k_r}_o(S_o)
= S_o^2$, rescaling $x\frac{\phi}{q}\rightarrow x$ and combining
(\ref{governing0})-(\ref{capillary}), the single equation for the
water saturation $u=S_w$ is
\begin{equation}
\frac{\partial u}{\partial t}+ \frac{\partial}{\partial
x}\left[\frac{u^2}{u^2+M(1-u)^2}\right]=-\frac{\partial}{\partial
x}\left[\frac{\phi^2}{q^2}\frac{k(1-u)^2u^2}{\mu_w(1-u)^2+\mu_ou^2}\frac{\partial}{\partial
x}\left(\frac{p_c(u)}{\phi}-\tau\frac{\partial u}{\partial
t}\right)\right] \label{governing1}
\end{equation}
where $M=\frac{\mu_w}{\mu_o}$ \cite{Ying}. Linearizing the right
hand side of (\ref{governing1}) and rescaling the equation as in
\cite{duijn, duijn2},
the modified Buckley-Leverett equation (MBL) is derived as
\begin{equation}
\frac{\partial u}{\partial t}+\frac{\partial f(u)}{\partial
x}=\epsilon\frac{\partial^2 u}{\partial
x^2}+\epsilon^2\tau\frac{\partial^3 u}{\partial x^2\partial t}
\label{MBL}
\end{equation}
where the water fractional flow function $f(u)$ is given as in
(\ref{flux}). Notice that, if $P_c$ in (\ref{capillary}) is taken to
be constant, then (\ref{governing1}) gives the classical BL
equation; while if the dispersive parameter $\tau$ is taken to be
zero, then (\ref{MBL}) gives the viscous BL equation, which still
displays monotone water saturation profile. Thus, in addition to
the classical second order viscous term $\epsilon u_{xx}$, the MBL
equation (\ref{MBL}) is an extension involving a third order mixed derivative term
$\epsilon^2\tau u_{xxt}$. Van Dujin et al. \cite{duijn} showed
that the value $\tau$ is critical in determining the type of the
solution profile. In particular, for certain Riemann problems, the solution
profile of (\ref{MBL}) is not monotone when $\tau$ is larger than
the threshold value $\tau_*$, where $\tau_*$ was numerically
 determined to be 0.61 \cite{duijn}. The non-monotonicity of the
solution profile is consistent with the experimental observations
\cite{Dicarlo}.

The classical BL equation (\ref{BL}) is hyperbolic, and the
numerical schemes for hyperbolic equations have been well
developed (e.g. \cite{Leveque, Leveque2, Cockburn2, Cockburn, NT, KL}
). The MBL
equation (\ref{MBL}), however, is pseudo-parabolic, we will
illustrate how to extend the central schemes 
\cite{NT,KL,KL2} to solve
(\ref{MBL}) numerically. Unlike the finite domain of dependence
for the classical BL equation (\ref{BL}), the domain of dependence for the MBL equation (\ref{MBL}) is infinite.
 This naturally raises the question for the choice of
computational domain. To answer this question, we will  first study the
MBL equation equipped with two types of domains and corresponding
boundary conditions. One is the half line boundary value problem
\begin{equation}
\begin{split}
u_t+(f(u))_x~ =~\epsilon u_{xx} + \epsilon^2\tau u_{xxt}
 \qquad&\mathrm{in}  \qquad Q=\{(x,t): x>0, t>0\} \\
u(x,0)~ =~ u_0(x) \qquad&\quad\qquad x\in[0,\infty)\\
u(0,t) ~=~g_u(t) ,\quad \lim_{x\rightarrow\infty}u(x,t)~=~0\qquad &\quad\qquad t\in[0,\infty)\\
u_0(0) ~=~g_u(0)
\qquad&\quad\qquad\mathrm{compatibility~condition}
\end{split}\label{quarter_plane_0}
\end{equation}
and the other one is finite interval boundary value problem
\begin{equation}
\begin{split}
v_t+(f(v))_x ~=~\epsilon v_{xx} + \epsilon^2\tau v_{xxt}
 \qquad&\mathrm{in}  \qquad \widetilde{Q}=\{(x,t): x\in(0,L), t>0\}\\
v(x,0) ~=~ v_0(x) \qquad&\quad\qquad  x\in[0,L]\\
v(0,t) ~=~g_v(t), \quad v(L,t) ~=~h(t)\qquad&\quad\qquad t\in[0,\infty)  \\
v_0(0) ~=~g_v(0),\quad v_0(L) ~=~ h(0)\qquad
&\quad\qquad\mathrm{compatibility~condition}. 
\end{split}\label{BV_0}
\end{equation}
Considering
\begin{equation}
u_0(x)=\left\{
\begin{array}{lll}
v_0(x) & \mathrm{for} & x\in[0,L]\\
0 &\mathrm{for} & x\in[L,+\infty)
\end{array}\right.
, \qquad g_u(t)=g_v(t)\equiv g(t), \qquad h(t)\equiv 0, \label{assumptions}
\end{equation}
we will show the relation between the solutions of problems
(\ref{quarter_plane_0}) and (\ref{BV_0}). To the best knowledge of
the authors, there is no such study for MBL equation (\ref{MBL}).
Similar questions were answered for BBM equation \cite{bona,
bona1}.

The organization of this paper is as follows. Section
\ref{comparison} will bring forward the exact theory comparing the
solutions of (\ref{quarter_plane_0}) and (\ref{BV_0}). The
difference between the solutions of these two types of problems
decays exponentially with respect to the length of the interval
$L$ for practically interesting initial profiles. This provides a
theoretical justification for the choice of the computational
domain. In section \ref{central}, high order central schemes will
be developed for MBL equation in finite interval domain. We
provide a detailed derivation on how to extend the central schemes
\cite{NT,KL} for conservation laws to solve the MBL equation
(\ref{MBL}). The idea of adopting numerical schemes originally
designed for hyperbolic equations to pseudo-parabolic equations is
not restricted to central type schemes only (\cite{Shu1,Shu2}).
The numerical results  in section \ref{results} show that the
water saturation profile strongly depends on the dispersive
parameter $\tau$ value as studied in \cite{duijn}. For
$\tau>\tau_*$, the MBL equation (\ref{MBL}) gives non-monotone
water saturation profiles for certain Riemann problems as
suggested by experimental observations \cite{Dicarlo}. Section
\ref{conclusion} gives the conclusion of the paper and the
possible future directions.

\section{The half line problem versus the finite interval problem}
\label{comparison} Let $u(x,t)$ be the solution to the half line
problem (\ref{quarter_plane_0}), and let $v(x,t)$ be the solution
to the finite interval problem (\ref{BV_0}). We
consider the natural assumptions (\ref{assumptions}). The goal of
this section is to develop an estimate of the difference between
$u$ and $v$ on the spatial interval $[0,L]$ at a given finite time
$t$.  The main result of this section is
\begin{my_thm}[The main Theorem]
If $u_0(x)$ satisfies
\begin{equation}
u_0(x)=\left\{
\begin{array}{llll}
 C_u & & &x\in [0,L_0]\\
0& & &x > L_0
\end{array}
\label{u0_new} \right.
\end{equation}
where $L_0<L$ and $C_u,$  
%
are positive constants, then
\begin{equation*}
\norm{ u(\cdot,t)-v(\cdot,t)}_{H_{L,\epsilon,\tau}^1}\le D_{1;\epsilon,\tau}(t)
e^{-\frac{\lambda L}{\epsilon\sqrt{\tau}}}+D_{2;\epsilon,\tau}(t)
e^{-\frac{\lambda (L-L_0)}{\epsilon\sqrt{\tau}}}
\end{equation*}
for some $0<\lambda<1$,
$D_{1;\epsilon,\tau}(t)>0$ and $D_{2;\epsilon,\tau}(t)>0$, where
\begin{equation*}
\norm{ Y(\cdot,t) }_{H_{L,\epsilon,\tau}^1}:=\sqrt{\int_0^L
Y(x,t)^2+(\epsilon\sqrt{\tau}Y_x(x,t))^2\, dx}
\end{equation*}
\label{main_new}
\end{my_thm}
Notice that
%
the initial condition (\ref{u0_new}) we considered is the Riemann
problem. Theorem \ref{main_new} shows that the solution to the
half line problem (\ref{quarter_plane_0}) can be approximated as
accurately as one wants by the solution to the finite interval
 problem (\ref{BV_0}) in the sense that
$D_{1;\epsilon,\tau}(t)$, $D_{2;\epsilon,\tau}(t)$, $\frac{\lambda
L}{\epsilon\sqrt{\tau}}$ and $\frac{\lambda
(L-L_0)}{\epsilon\sqrt{\tau}}$ can be controlled.

To prove theorem \ref{main_new}, we first derive the implicit solution formulae for the half line problem and the finite interval problem in section \ref{qp} and section \ref{tp} respectively. The implicit solution formulae are in integral form, which are derived by separating the $x$-derivative from the $t$-derivative, and formally solving a first order linear ODE in $t$ and a second order non-homogeneous ODE in $x$. In section \ref{comp}, we use Gronwall's  inequality multiple times to obtain the desired result in theorem \ref{main_new}.

\subsection{Half line problem}
\label{qp} In this section, we derive the implicit solution
formula for the half line problem (\ref{quarter_plane_0}) (with $g_u(t) = g(t)$). 
To solve (\ref{quarter_plane_0}), we first  rewrite (\ref{quarter_plane_0}) by separating the $x$-derivative from the
$t$-derivative,
\begin{equation}
\left(I -\epsilon^2\tau\frac{\partial^2}{\partial
x^2}\right)\left(u_t+\frac{1}{\epsilon\tau}u\right) = \frac{1}{\epsilon\tau} u-(f(u))_x 
 .\label{sep}
\end{equation}
By using integrating factor method, we formally  integrate (\ref{sep}) over $[0,t]$ to obtain
\begin{equation}
\left(I -\epsilon^2\tau\frac{\partial^2}{\partial
x^2}\right)\left(u - e^{-\frac{t}{\epsilon\tau}}u_0\right)
=\int_0^t \left(\frac{1}{\epsilon\tau}
u -(f(u))_x \right)e^{-\frac{t-s}{\epsilon\tau}}\, ds .\label{t_done}
\end{equation}
Furthermore, we let
\begin{equation}
A = u - e^{-\frac{t}{\epsilon\tau}}u_0 ,\label{eq1prime}
\end{equation}
then (\ref{t_done}) can be written as
\begin{equation}
 A'' -\frac{1}{ \epsilon^2\tau}A = \int_0^t \left(-\frac{1}{\epsilon^3\tau^2} u+\frac{1}{\epsilon^2\tau}(f(u))_x 
\right) e^{-\frac{t-s}{\epsilon\tau}}\,ds,
\qquad \mathrm{where} ~~ '=\frac{\partial}{\partial x}.
\label{eq2}
\end{equation}
Notice that  (\ref{eq2}) is a second-order  non-homogeneous ODE in
$x$-variable along with the boundary conditions
\begin{equation}
\begin{split}
A(0,t) &= u(0,t)-e^{-\frac{t}{\epsilon\tau}}u_0(0) ~=~g(t)-e^{-\frac{t}{\epsilon\tau}}g(0) ,\\
A(\infty,t)&=u(\infty,t)-e^{-\frac{t}{\epsilon\tau}}u_0(\infty)~=~0.
\end{split} \label{BCs}
\end{equation}
To solve (\ref{eq2}), we first solve the corresponding linear
homogeneous equation with the non-zero boundary conditions
(\ref{BCs}). We then find a particular solution for the
non-homogeneous equation with zero boundary conditions by
introducing a Green's function $G(x,\xi)$ and a kernel $K(x,\xi)$
for the non-homogeneous terms $
u$ and
$(f(u))_x$ respectively. Combining the solutions for the two
non-homogeneous terms and the homogeneous part with boundary
conditions, we get the solution for equation (\ref{eq2})
satisfying the boundary conditions (\ref{BCs}):
\begin{equation}
\begin{split}
A(x,t) &=  -\frac{1}{\epsilon^3\tau^2}\int_0^t\int_0^{+\infty}G(x,\xi)u(\xi,s)e^{-\frac{t-s}{\epsilon\tau}}\,d\xi\,ds  \\
&\quad +\frac{1}{\epsilon^2\tau}\int_0^t\int_0^{+\infty}
K(x,\xi)f(u)e^{-\frac{t-s}{\epsilon\tau}}\,d\xi\,ds \\
&\quad + \left(g(t)-e^{-\frac{t}{\epsilon\tau}}g(0)\right)e^{-\frac{x}{\epsilon\sqrt{\tau}}}\end{split}
\label{quarter_plane_sol}
\end{equation}
where the Green's function $G(x,\xi)$ and the kernel $K(x,\xi)$
are
\begin{eqnarray}
G(x,\xi) &=&\frac{\epsilon\sqrt{\tau}}{2}
\left(e^{-\frac{x+\xi}{\epsilon\sqrt{\tau}}}-e^{-\frac{|x-\xi|}{\epsilon\sqrt{\tau}}}\right),\label{G1}\\
K(x,\xi) &=&-\frac{\partial G(x,\xi)}{\partial\xi}~=~
\frac{1}{2}\left(e^{-\frac{x+\xi}{\epsilon\sqrt{\tau}}}+\mathrm{sgn}(x-\xi)e^{-\frac{|x-\xi|}{\epsilon\sqrt{\tau}}}\right).\label{K2}
\end{eqnarray}
To recover the solution for the half line problem (\ref{quarter_plane_0}), we refer to the definition of $A$ in (\ref{eq1prime}). Thus, the implicit solution formula for the half line problem (\ref{quarter_plane_0}) is
\begin{equation}
\begin{split}
u(x,t) &=-\frac{1}{2\epsilon^2\tau\sqrt{\tau}}\int_0^t\int_0^{+\infty}\left(e^{-\frac{x+\xi}{\epsilon\sqrt{\tau}}}-e^{-\frac{|x-\xi|}{\epsilon\sqrt{\tau}}}\right)u(\xi,s)e^{-\frac{t-s}{\epsilon\tau}}\,d\xi\,ds
 \\
&\quad+\frac{1}{2\epsilon^2\tau}\int_0^t\int_0^{+\infty}\left(e^{-\frac{x+\xi}{\epsilon\sqrt{\tau}}}+\mathrm{sgn}(x-\xi)e^{-\frac{|x-\xi|}{\epsilon\sqrt{\tau}}}\right)f(u)e^{-\frac{t-s}{\epsilon\tau}}\,d\xi\,ds\\
&\quad+\left(g(t)-e^{-\frac{t}{\epsilon\tau}}g(0)\right)e^{-\frac{x}{\epsilon\sqrt{\tau}}}
+ e^{-\frac{t}{\epsilon\tau}}u_0(x)
.
\end{split}\label{quarter_plane_sol_final}
\end{equation}

\subsection{Finite interval problem}
\label{tp}
The implicit solution for the finite interval problem (\ref{BV_0}) (with $g_v(t) = g(t)$) 
can be solved in a similar way.
The only difference is that the additional boundary condition $h(t)$ at $x=L$ in (\ref{BV_0}) gives different boundary conditions for the non-homogeneous ODE in $x$-variable. Denote
\begin{equation}
A^L = v - e^{-\frac{t}{\epsilon\tau}}v_0 ,\label{BVeq1prime}
\end{equation}
then it satisfies
\begin{equation}
(A^L)''-\frac{1}{\epsilon^2\tau}A^L =\int_0^t \left(-\frac{1}{\epsilon^3\tau^2} v+\frac{1}{\epsilon^2\tau}(f(v))_x 
\right) e^{-\frac{t-s}{\epsilon\tau}}\,ds
\quad \mathrm{where} \quad'=\frac{\partial}{\partial x}\label{BVeq2}
\end{equation}
with the boundary conditions
\begin{equation*}
\begin{split}
A^L(0,t) &= v(0,t)-e^{-\frac{t}{\epsilon\tau}}v_0(0) =g(t)-e^{-\frac{t}{\epsilon\tau}}g(0), \\
A^L(L,t)&=v(L,t)-e^{-\frac{t}{\epsilon\tau}}v_0(L)
=h(t)-e^{-\frac{t}{\epsilon\tau}}h(0).
\end{split}
\end{equation*}

These boundary conditions affect both the homogeneous solution and
the particular solution  of (\ref{BVeq2}) as follows
\begin{equation}
\begin{split}
A^L(x,t)&=-\frac{1}{\epsilon^3\tau^2}\int_0^t\int_0^{L}G^L(x,\xi)v(\xi,s)e^{-\frac{t-s}{\epsilon\tau}}\,d\xi\,ds \\
&\quad+\frac{1}{\epsilon^2\tau}\int_0^t\int_0^{L}K^L(x,\xi)f(v)e^{-\frac{t-s}{\epsilon\tau}}\,d\xi\,ds \\
&\quad+ c_1(t) \phi_1(x)+c_2(t) \phi_2(x)
\end{split}\label{BV_sol}
\end{equation}
where the Green's function $G^L(x,\xi)$, the kernel $K^L(x,\xi)$
and the bases for the homogeneous solutions are
\begin{equation}
G^L(x,\xi)=\frac{\epsilon\sqrt{\tau}}{2(e^{\frac{2L}{\epsilon\sqrt{\tau}}}-1)}
\left(e^{\frac{x+\xi}{\epsilon\sqrt{\tau}}}+
e^{\frac{2L-(x+\xi)}{\epsilon\sqrt{\tau}}}
-e^{\frac{|x-\xi|}{\epsilon\sqrt{\tau}}}
-e^{\frac{2L-|x-\xi|}{\epsilon\sqrt{\tau}}} \right), \label{G1L}
\end{equation}
\begin{equation}
\begin{split}
K^L(x,\xi) &= -\frac{1}{2(e^{\frac{2L}{\epsilon\sqrt{\tau}}}-1)}\left(e^{\frac{x+\xi}{\epsilon\sqrt{\tau}}}
-e^{\frac{2L-(x+\xi)}{\epsilon\sqrt{\tau}}} \right.\\
&\qquad\qquad\qquad\left. +\mathrm{sgn}(x-\xi)e^{\frac{|x-\xi|}{\epsilon\sqrt{\tau}}}
-\mathrm{sgn}(x-\xi)e^{\frac{2L-|x-\xi|}{\epsilon\sqrt{\tau}}}\right),
\end{split}
\label{K2L}
\end{equation}
\begin{align}
&c_1(t)=g(t)-e^{-\frac{t}{\epsilon\tau}}g(0) , \qquad\qquad\quad c_2(t)=h(t)-e^{-\frac{t}{\epsilon\tau}}h(0),\label{cs}\\
 &\phi_1(x) =
\frac{e^{\frac{L-x}{\epsilon\sqrt{\tau}}}-e^{\frac{-L+x}{\epsilon\sqrt{\tau}}}}
{e^{\frac{L}{\epsilon\sqrt{\tau}}}-e^{-\frac{L}{\epsilon\sqrt{\tau}}}},
\qquad\text{and}\quad \phi_2(x) =
\frac{e^{\frac{x}{\epsilon\sqrt{\tau}}}-e^{-\frac{x}{\epsilon\sqrt{\tau}}}}
{e^{\frac{L}{\epsilon\sqrt{\tau}}}-e^{-\frac{L}{\epsilon\sqrt{\tau}}}}.
\label{phi}
\end{align}
Thus, the implicit solution formula for the finite interval problem (\ref{BV_0}) is
\begin{equation}
\begin{split}
v(x,t)=&
-\frac{1}{2\epsilon^2\tau\sqrt{\tau}(e^{\frac{2L}{\epsilon\sqrt{\tau}}}-1)}
\int_0^t\int_0^{L}\left(e^{\frac{x+\xi}{\epsilon\sqrt{\tau}}}+
e^{\frac{2L-(x+\xi)}{\epsilon\sqrt{\tau}}}
-e^{\frac{|x-\xi|}{\epsilon\sqrt{\tau}}} \right.\\
&\hspace{5cm} 
\left. -e^{\frac{2L-|x-\xi|}{\epsilon\sqrt{\tau}}}
\right) v(\xi,s)e^{-\frac{t-s}{\epsilon\tau}}\,d\xi\,ds\\
&-\frac{1}{2\epsilon^2\tau(e^{\frac{2L}{\epsilon\sqrt{\tau}}}-1)}\int_0^t\int_0^{L}\left(e^{\frac{x+\xi}{\epsilon\sqrt{\tau}}}
-e^{\frac{2L-(x+\xi)}{\epsilon\sqrt{\tau}}}
+\mathrm{sgn}(x-\xi)e^{\frac{|x-\xi|}{\epsilon\sqrt{\tau}}}\right.\\
&\left.\hspace{4cm}
-\mathrm{sgn}(x-\xi)e^{\frac{2L-|x-\xi|}{\epsilon\sqrt{\tau}}}
\right)
f(v)e^{-\frac{t-s}{\epsilon\tau}}\,d\xi\,ds\\
&+
c_1(t)\phi_1(x)+
c_2(t)\phi_2(x)
 + e^{-\frac{t}{\epsilon\tau}}v_0(x) . 
\end{split}
\label{BV_sol_final}
\end{equation}

\subsection{Comparisons}
\label{comp} 
In this section, we will prove 
that the solution  $u(x,t)$ to the half line problem can
be approximated as accurately as one wants by the solution
$v(x,t)$ to the finite interval problem as stated
in Theorem \ref{main_new}.
%

Due to the difference in the integration domains, we do not use
\eqref{quarter_plane_sol_final} and \eqref{BV_sol_final} directly for the comparison.
Instead, we decompose
$u(x,t)$ ($v(x,t)$ respectively) into two parts: $U(x,t)$ and $u_L(x,t)$
 ($V(x,t)$ and $v_L(x,t)$ respectively), 
such that $U(x,t)$ ($V(x,t)$ respectively) enjoys zero
initial condition and boundary conditions at $x=0$ and $x=L$.
We estimate the difference between $u(\cdot,t)$ and $v(\cdot,t)$ by estimating the differences between
$u_L(\cdot,t)$ and $v_L(\cdot,t)$, $U(\cdot,t)$ and $V(\cdot,t)$,   then applying  the triangle inequality.

\textcolor{red}{
}

\subsubsection{Definitions and lemmas}
\label{def}
To assist the proof of Theorem \ref{main_new} in section \ref{proof}, we introduce some new notations in this section. We first decompose $u(x,t)$ as sum of two terms $U(x,t)$ and $u_L(x,t)$, such that
\begin{equation*}
u(x,t)=U(x,t)+u_L(x,t)\qquad x\in[0,+\infty)
\end{equation*}
where
\begin{equation}
u_L = e^{-\frac{t}{\epsilon\tau}}u_0(x) + c_1(t)
e^{-\frac{x}{\epsilon\sqrt{\tau}}} +
\left(u(L,t)-c_1(t)e^{-\frac{L}{\epsilon\sqrt{\tau}}}-e^{-\frac{t}{\epsilon\tau}}u_0(L)\right)\phi_2(x)
\label{uL}
\end{equation}
and $ c_1(t)$ and $\phi_2(x)$ are given in (\ref{cs}) and
(\ref{phi}) respectively. With this definition, $u_L$ takes care
of the initial condition $u_0(x)$ and boundary conditions $g(t)$
at $x=0$ and $x=L$ for $u(x,t)$.
Then $U$ satisfies an equation slightly different from the
equation $u$ satisfies in (\ref{quarter_plane_0}):
\begin{equation}
\begin{split}
U_t-\epsilon U_{xx}-\epsilon^2\tau U_{xxt}& =\left(u_t-\epsilon u_{xx}-\epsilon^2\tau u_{xxt}\right) -\left((u_L)_t-\epsilon (u_L)_{xx}-\epsilon^2\tau (u_L)_{xxt}\right) \\
&=-\left(f(u)\right)_x+\frac{1}{\epsilon\tau}u_L(x,t) 
\end{split}\label{Ueq}
%
\end{equation}
In addition, $U(x,t)$ has zero initial condition and boundary conditions at $x=0$ and $x=L$, i.e.,
\begin{eqnarray}
U(x,0) =0,  \qquad U(0,t) =0, \qquad U(L,t) =0.
\label{Uzero}
\end{eqnarray}
Similarly, for $v(x,t)$, let
\begin{equation*}
v(x,t) =  V(x,t) + v_L(x,t) \qquad x\in[0,L]
\end{equation*}
where
\begin{equation}
v_L = e^{-\frac{t}{\epsilon\tau}}v_0(x) + c_1(t) \phi_1(x) +
c_2(t)\phi_2(x) \label{vL}
\end{equation}
and $c_1(t)$, $c_2(t)$ and $\phi_1(x)$, $\phi_2(x)$ are given in
(\ref{cs}) and (\ref{phi}) respectively. With this definition,
$v_L$ takes care of the initial condition $v_0(x)$ and boundary
conditions $g(t)$ and $h(t)$ at $x=0$ and $x=L$ for $v(x,t)$.
Then $V$ satisfies an equation slightly different from the
equation $v$ satisfies in (\ref{BV_0}):
\begin{equation}
\begin{split}
V_t-\epsilon V_{xx}-\epsilon^2\tau V_{xxt}
&=-\left(f(v)\right)_x+\frac{1}{\epsilon\tau}v_L(x,t) 
\end{split}\label{Veq}
%
\end{equation}
with
\begin{eqnarray}
V(x,0) =0,  \qquad V(0,t) =0, \qquad V(L,t) =0.
\label{Vzero}
\end{eqnarray}
Since, in the end, we want to study the difference between $U(x,t)$ and $V(x,t)$, we define
\begin{equation*}
W(x,t) = V(x,t) - U(x,t) \qquad\mathrm{for}\qquad x\in[0,L].
\end{equation*}
Because of (\ref{Ueq}) and (\ref{Veq}), we have
\begin{equation}
W_t-\epsilon W_{xx}-\epsilon^2\tau W_{xxt} =
-\left(f(v)-f(u)\right)_x +\frac{1}{\epsilon\tau}(v_L-u_L) .
\label{Weq}
\end{equation}
In lieu of (\ref{Uzero}) and (\ref{Vzero}), $W(x,t)$ also has zero initial condition and boundary conditions at $x=0$ and $x=L$, i.e.,
\begin{eqnarray}
W(x,0) =0,  \qquad W(0,t) =0, \qquad W(L,t) =0. \label{W_BC}
\end{eqnarray}

Now, to estimate $\norm{u-v}$, we can estimate
$\norm{W}=\norm{V-U}$  and estimate $\norm{u_L-v_L}$ separately.
These estimates are done in section \ref{proof}.

Next, we state the lemmas needed in the proof of Theorem
\ref{main_new}. The proof of the lemmas can be found in the
appendix \ref{proof_lemma} and \cite{Ying}. In all the lemmas, we
assume $0<\lambda<1$ and $u_0(x)$ satisfies
\begin{eqnarray}
u_0(x)=\left\{
\begin{array}{llll}
 C_u & & &x\in [0,L_0]\\
0& & &x> L_0
\end{array}
\label{u0} \right.
\end{eqnarray}
where $L_0<L$ and $C_u$ 
are positive constants. Notice that the constraint
$\lambda\in(0,1)$ is crucial in Lemmas \ref{lemma_old1},
\ref{lemma_old2}.

\begin{my_lem}
$f(u) = \frac{u^2}{u^2+M(1-u)^2} \le Du $ ~where~
$D=\frac{f(\alpha)}{\alpha}$ and $\alpha=\sqrt{\frac{M}{M+1}}$. \label{lemmafu}
\end{my_lem}

\begin{my_lem}
\label{lemma_old1}
\begin{enumerate}[(i)]
\item\label{lemma1} $\int_{0}^{+\infty}
\left|e^{-\frac{x+\xi}{\epsilon\sqrt{\tau}}}-e^{-\frac{|x-\xi|}{\epsilon\sqrt{\tau}}}\right|
e^{\frac{\lambda x -\lambda\xi}{\epsilon\sqrt{\tau}}}\,d\xi \le
\frac{2\epsilon\sqrt{\tau}}{1-\lambda^2} . $ 
\item \label{lemma2}
$\int_{0}^{+\infty}
\left|e^{-\frac{x+\xi}{\epsilon\sqrt{\tau}}}-e^{-\frac{|x-\xi|}{\epsilon\sqrt{\tau}}}\right|
e^{\frac{\lambda x-\xi}{\epsilon\sqrt{\tau}}}\,d\xi \le
\frac{\epsilon\sqrt{\tau}}{e(1-\lambda)}$ . 
\item\label{lemma3}
$\int_{0}^{+\infty}\left|e^{-\frac{x+\xi}{\epsilon\sqrt{\tau}}}-e^{-\frac{|x-\xi|}{\epsilon\sqrt{\tau}}}\right|e^{\frac{\lambda
x}{\epsilon\sqrt{\tau}}} |u_0(\xi)|\,d\xi \le
2C_u\epsilon\sqrt{\tau}e^{\frac{\lambda
L_0}{\epsilon\sqrt{\tau}}}
$
.
\end{enumerate}
\end{my_lem}

\begin{my_lem}
\label{lemma_old2}
\begin{enumerate}[(i)]
\item\label{lemma11} $\int_{0}^{+\infty}
\left|e^{-\frac{x+\xi}{\epsilon\sqrt{\tau}}}+\mathrm{sgn}(x-\xi)e^{-\frac{|x-\xi|}{\epsilon\sqrt{\tau}}}\right|
e^{\frac{\lambda
x-\lambda\xi}{\epsilon\sqrt{\tau}}}\,d\xi\le\frac{2\epsilon\sqrt{\tau}}{1-\lambda^2}$
. \item \label{lemma12} $\int_{0}^{+\infty} \left|
e^{-\frac{x+\xi}{\epsilon\sqrt{\tau}}}+\mathrm{sgn}(x-\xi)e^{-\frac{|x-\xi|}{\epsilon\sqrt{\tau}}}\right|
e^{\frac{\lambda x-\xi}{\epsilon\sqrt{\tau}}}\,d\xi\le
\epsilon\sqrt{\tau}+\frac{\epsilon\sqrt{\tau}}{e(1-\lambda)}$ .
\item\label{lemma13} $\int_{0}^{+\infty}
\left|e^{-\frac{x+\xi}{\epsilon\sqrt{\tau}}}+\mathrm{sgn}(x-\xi)e^{-\frac{|x-\xi|}{\epsilon\sqrt{\tau}}}\right|
e^{\frac{\lambda x}{\epsilon\sqrt{\tau}}}|u_0(\xi)|\,d\xi\le
2C_u\epsilon\sqrt{\tau}e^{\frac{\lambda L_0}{\epsilon\sqrt{\tau}}}
$ .
\end{enumerate}
\end{my_lem}

\begin{my_lem}
\label{lemma_old4}
\begin{enumerate}[(i)]
\item\label{lemmaphi1}
$\left|\phi_1(x)-e^{-\frac{x}{\epsilon\sqrt{\tau}}}\right|=
e^{-\frac{L}{\epsilon\sqrt{\tau}}}\left|\phi_2(x)\right|$ .
\item\label{lemmaphi2}
$\left|\phi_2(x)\right| \le 1$ ~for~  $x\in[0,L]$ .
\item\label{lemmaphi2prime}
$\left|\phi_2'(x)\right|\le \frac{2}{\epsilon\sqrt{\tau}}$ ~if~ $\epsilon\ll 1$ ~for~ $x\in[0,L]$ .
\end{enumerate}
\end{my_lem}

\noindent Last but not least, the norm that we will use in Theorem \ref{main_new} and its proof is
\begin{equation}
\norm{ Y(\cdot,t) }_{H_{L,\epsilon,\tau}^1}:=\sqrt{\int_0^L
Y(x,t)^2+(\epsilon\sqrt{\tau}Y_x(x,t))^2\, dx} .\label{Ynorm}
\end{equation}

\subsubsection{A proposition}
\label{prop}
In this section, we will give a critical estimate, which is essential in the calculation of maximum difference $\norm{u_L(\cdot,t)-v_L(\cdot,t)}_\infty$ in section \ref{proof}.
By comparing $u_L(x,t)$ and $v_L(x,t)$ given in (\ref{uL}) and (\ref{vL}) respectively, it is clear that the coefficient $u(L,t)-c_1(t)e^{-\frac{L}{\epsilon\sqrt{\tau}}}-e^{-\frac{t}{\epsilon\tau}}u_0(L)$ for $\phi_2(x)$ appeared in (\ref{uL}) needs to be compared with the corresponding coefficient $c_2(t)$ for $\phi_2(x)$ appeared in (\ref{vL}).
We thus  define a space-dependent function
\begin{equation}
U_{c_2}(x,t) = u(x,t) -
c_1(t)e^{-\frac{x}{\epsilon\sqrt{\tau}}}-
e^{-\frac{t}{\epsilon\tau}}u_0(x)\label{Uc2def}
\end{equation}
and establish the following proposition
\begin{my_prop}
\begin{equation}
\left|U_{c_2}(L,t)\right|
\le
a_{\tau}(t)e^\frac{b_{\tau}t}{\epsilon\tau}
e^{-\frac{\lambda L}{\epsilon\sqrt{\tau}}} +
c_{\tau}\frac{t}{\epsilon\tau}e^\frac{(b_{\tau}-1)t}{\epsilon\tau}e^{-\frac{\lambda(L-L_0)}{\epsilon\sqrt{\tau}}}\label{prop_est}
\end{equation}
for some parameter-dependent constants $a_{\tau}$, $b_{\tau}$ and $c_{\tau}$.
\end{my_prop}
\begin{proof}
Based on the implicit solution formula
(\ref{quarter_plane_sol_final}) derived in section \ref{qp}, Lemma
\ref{lemmafu} and the relationship between $U_{c_2}$ and $u$ given
in (\ref{Uc2def}), we can get an inequality in terms of $U_{c_2}$
\begin{equation}
\begin{split}
&\left|U_{c_2}(x,t)\right|  \le  \frac{1}{2\epsilon^2\tau\sqrt{\tau}}
\left[\int_0^t\int_0^{+\infty}
\left|e^{-\frac{x+\xi}{\epsilon\sqrt{\tau}}}-e^{-\frac{|x-\xi|}{\epsilon\sqrt{\tau}}}\right|
\left|U_{c_2}(\xi,s)\right|e^{-\frac{t-s}{\epsilon\tau}}\,d\xi\,ds
\right. \\
&\qquad+\int_0^t\int_0^{+\infty}\left|e^{-\frac{x+\xi}{\epsilon\sqrt{\tau}}}-e^{-\frac{|x-\xi|}{\epsilon\sqrt{\tau}}}\right|\left|c_1(s)\right|e^{-\frac{\xi}{\epsilon\sqrt{\tau}}}e^{-\frac{t-s}{\epsilon\tau}}\,d\xi\,ds
 \\
&\qquad+\left.\int_0^t\int_0^{+\infty}\left|e^{-\frac{x+\xi}{\epsilon\sqrt{\tau}}}-e^{-\frac{|x-\xi|}{\epsilon\sqrt{\tau}}}\right|\left|u_0(\xi)\right|e^{-\frac{t}{\epsilon\tau}}\,d\xi\,ds
\right] \\
&\quad+\frac{D}{2\epsilon^2\tau}\left[\int_0^t\int_0^{+\infty}\left|e^{-\frac{x+\xi}{\epsilon\sqrt{\tau}}}+\mathrm{sgn}(x-\xi)e^{-\frac{|x-\xi|}{\epsilon\sqrt{\tau}}}\right|\left|U_{c_2}(\xi,s)\right|e^{-\frac{t-s}{\epsilon\tau}}\,d\xi\,ds\right. \\
&\qquad+\int_0^t\int_0^{+\infty}\left|e^{-\frac{x+\xi}{\epsilon\sqrt{\tau}}}+\mathrm{sgn}(x-\xi)e^{-\frac{|x-\xi|}{\epsilon\sqrt{\tau}}}\right|\left|c_1(s)\right|e^{-\frac{\xi}{\epsilon\sqrt{\tau}}}e^{-\frac{t-s}{\epsilon\tau}}\,d\xi\,ds \\
&\qquad+\left.\int_0^t\int_0^{+\infty}\left|e^{-\frac{x+\xi}{\epsilon\sqrt{\tau}}}+\mathrm{sgn}(x-\xi)e^{-\frac{|x-\xi|}{\epsilon\sqrt{\tau}}}\right|\left|u_0(\xi)\right|e^{-\frac{t}{\epsilon\tau}}\,d\xi\,ds\right] .
\end{split}\label{Uc2ineq}
\end{equation}
To show that $U_{c_2}(x,t)$ decays exponentially with respect to
$x$, we pull out an exponential term by writing
$U_{c_2}(x,t) = e^{-\frac{\lambda
x}{\epsilon\sqrt{\tau}}}e^{-\frac{t}{\epsilon\tau}}\tilde{U}(x,t)$, where $0<\lambda<1$, such
that
\begin{equation}
\tilde{U}(x,t) = e^{\frac{\lambda x}{\epsilon\sqrt{\tau}}}e^{\frac{t}{\epsilon\tau}}U_{c_2}(x,t) ,
\label{utilde}
\end{equation}
then (\ref{Uc2ineq}) can be rewritten in terms of $\tilde{U}(x,t)$
as follows
\begin{equation}
\begin{split}
&\left|\tilde{U}(x,t)\right| \le 
\frac{1}{2\epsilon^2\tau\sqrt{\tau}}
\left[\int_0^t\int_0^{+\infty}
\left|e^{-\frac{x+\xi}{\epsilon\sqrt{\tau}}}-e^{-\frac{|x-\xi|}{\epsilon\sqrt{\tau}}}\right|
e^{\frac{\lambda
x-\lambda\xi}{\epsilon\sqrt{\tau}}}\left|\tilde{U}(\xi,s)\right|
\,d\xi\,ds
\right.  \\
&\qquad+\int_0^t\int_0^{+\infty}\left|e^{-\frac{x+\xi}{\epsilon\sqrt{\tau}}}-e^{-\frac{|x-\xi|}{\epsilon\sqrt{\tau}}}\right|\left|c_1(s)\right|e^{\frac{\lambda
x-\xi}{\epsilon\sqrt{\tau}}}e^{\frac{s}{\epsilon\tau}}\,d\xi\,ds
 \\
&\qquad+\left.\int_0^t\int_0^{+\infty}\left|e^{-\frac{x+\xi}{\epsilon\sqrt{\tau}}}-e^{-\frac{|x-\xi|}{\epsilon\sqrt{\tau}}}\right|e^{\frac{\lambda
x}{\epsilon\sqrt{\tau}}}\left|u_0(\xi)\right|
\,d\xi\,ds
\right]  \\
&+\frac{D}{2\epsilon^2\tau}\left[\int_0^t\int_0^{+\infty}\left|e^{-\frac{x+\xi}{\epsilon\sqrt{\tau}}}+\mathrm{sgn}(x-\xi)e^{-\frac{|x-\xi|}{\epsilon\sqrt{\tau}}}\right|e^{\frac{\lambda x -\lambda\xi}{\epsilon\sqrt{\tau}}}\left|\tilde{U}_(\xi,s)\right|
\,d\xi\,ds\right.  \\
&\qquad+\int_0^t\int_0^{+\infty}\left|e^{-\frac{x+\xi}{\epsilon\sqrt{\tau}}}+\mathrm{sgn}(x-\xi)e^{-\frac{|x-\xi|}{\epsilon\sqrt{\tau}}}\right|\left|c_1(s)\right|e^{\frac{\lambda x-\xi}{\epsilon\sqrt{\tau}}}e^{\frac{s}{\epsilon\tau}}\,d\xi\,ds \\
&\qquad+\left.\int_0^t\int_0^{+\infty}\left|e^{-\frac{x+\xi}{\epsilon\sqrt{\tau}}}+\mathrm{sgn}(x-\xi)e^{-\frac{|x-\xi|}{\epsilon\sqrt{\tau}}}\right|e^{\frac{\lambda
x}{\epsilon\sqrt{\tau}}}\left|u_0(\xi)\right|
\,d\xi\,ds\right] .
\end{split}
\label{Utildeineq}
\end{equation}
Because of Lemmas \ref{lemma_old1}--\ref{lemma_old2}, we can get
the following estimate for
$\left|\tilde{U}(\cdot,t)\right|_{\infty}$ based on
(\ref{Utildeineq}) :
\begin{equation}
\begin{split}
&\left|\tilde{U}(\cdot,t)\right|_{\infty}\le
\frac{1}{2\epsilon^2\tau\sqrt{\tau}}\left[\frac{2\epsilon\sqrt{\tau}}{1-\lambda^2}\int_{0}^{t}|\tilde{U}(\cdot,s)|_\infty\,ds
 +\frac{\epsilon\sqrt{\tau}}{e(1-\lambda)}\int_{0}^{t}|c_1(s)|e^\frac{s}{\epsilon\tau}\,
 ds \right.  \\
&\qquad\qquad \left. +
2C_u\epsilon\sqrt{\tau}e^{\frac{\lambda L_0}{\epsilon\sqrt{\tau}}}
\int_{0}^{t} 1 
\, ds\right]\ \\
&\quad +\frac{D}{2\epsilon^2\tau}
\left[\frac{2\epsilon\sqrt{\tau}}{1-\lambda^2}\int_{0}^{t}
|\tilde{U}(\cdot,s)|_\infty \, ds +
\epsilon\sqrt{\tau}\left(1+\frac{1}{e(1-\lambda)}\right)\int_{0}^{t}|c_1(s)|e^\frac{s}{\epsilon\tau}
\, ds  \right. \\
&\qquad \qquad \left. +
2C_u\epsilon\sqrt{\tau}e^{\frac{\lambda
L_0}{\epsilon\sqrt{\tau}}}
\int_{0}^{t}1
\, ds
\right] \\
&\le\int_{0}^{t}\frac{b_{\tau}}{\epsilon\tau} |\tilde{U}(\cdot,s)|_\infty\,ds +
\int_{0}^{t} \frac{\tilde{a}_{\tau}(s)}{\epsilon\tau}\,ds
\end{split}
\label{Utildeest}
\end{equation}
where
\begin{equation*}
\begin{split}
b_{\tau}&=\frac{1+D\sqrt{\tau}}{1-\lambda^2},\qquad\qquad
\tilde{a}_{\tau}(t)
=a_{\tau}e^\frac{t}{\epsilon\tau}+c_{\tau} e^{\frac{\lambda L_0}{\epsilon\sqrt{\tau}}},\\
a_{\tau}&=\frac{|c_1(\cdot)|_\infty(1+D\sqrt{\tau}(e(1-\lambda)+1))}{2 e (1-\lambda)}
,\quad c_{\tau}=C_u(1+D\sqrt{\tau}).
\end{split}
\end{equation*}
By Gronwall's inequality, inequality (\ref{Utildeest})  gives that 
\begin{equation*}
\left|\tilde{U}(\cdot,t)\right|_{\infty}\le \int_{0}^{t}
\frac{\tilde{a}_{\tau}(t-s)}{\epsilon\tau} e^\frac{b_{\tau}(t-s)}{\epsilon\tau}\, ds
~\le ~\left(a_\tau e^\frac{t}{\epsilon\tau} + c_{\tau}
\frac{t}{\epsilon\tau} e^{\frac{\lambda
L_0}{\epsilon\sqrt{\tau}}}\right) e^\frac{b_{\tau}t}{\epsilon\tau}
\end{equation*}
Hence $\left|U_{c_2}(x,t)\right| \le
\left|\tilde{U}(\cdot,t)\right|_{\infty} e^{\frac{-\lambda
x}{\epsilon\sqrt{\tau}}}e^{-\frac{t}{\epsilon\tau}}
\le \left(a_\tau e^\frac{t}{\epsilon\tau} +c_{\tau} 
\frac{t}{\epsilon\tau} e^{\frac{\lambda
L_0}{\epsilon\sqrt{\tau}}}\right) e^\frac{b_{\tau}t}{\epsilon\tau}  e^{\frac{-\lambda
x}{\epsilon\sqrt{\tau}}}e^{-\frac{t}{\epsilon\tau}}$
i.e., $U_{c_2}(x,t)$ decays
exponentially with respect to $x$. In particular, when $x=L$, we have
\begin{equation}
 \left|U_{c_2}(L,t)\right| \le 
a_\tau e^\frac{b_\tau t}{\epsilon\tau}e^{-\frac{\lambda L}{\epsilon\sqrt{\tau}}}+c_\tau\frac{t}{\epsilon\tau}e^\frac{(b_\tau-1)t}{\epsilon\tau}e^{-\frac{\lambda(L-L_0)}{\epsilon\sqrt{\tau}}}
\label{Uc2}
\end{equation}
as given in (\ref{prop_est}).
\end{proof}

\subsubsection{Proof of Theorem \ref{main_new}}
\label{proof}
In this section, we will first find the maximum difference of 
$\norm{u_L(\cdot,t)-v_L(\cdot,t)}_\infty$, then we will derive 
$\norm{u_L(\cdot,t)-v_L(\cdot,t)}_{H_{L,\epsilon,\tau}^1}$ 
and $\norm{W(\cdot,t)}_{H_{L,\epsilon,\tau}^1}=\norm{U(\cdot,t)-V(\cdot,t)}_{H_{L,\epsilon,\tau}^1}$.
Combining these two, we will get an estimate for $\norm{u(\cdot,t)-v(\cdot,t)}_{H_{L,\epsilon,\tau}^1}$.
\begin{my_prop}
If $u_0(x)$ satisfies (\ref{u0}), then 
\begin{equation*}
\norm{u_L-v_L}_\infty\le
E_{1;\epsilon,\tau}(t)e^{-\frac{\lambda
L}{\epsilon\sqrt{\tau}}}+E_{2;\epsilon,\tau}(t)e^{-\frac{\lambda(L-L_0)}{\epsilon\sqrt{\tau}}}
\end{equation*}
 where $E_{1;\epsilon,\tau}(t)=
|c_1(\cdot)|_\infty+a_\tau e^\frac{b_\tau t}{\epsilon\tau}$ and $
E_{2;\epsilon,\tau}(t)=c_{\tau}\frac{t}{\epsilon\tau}e^\frac{(b_{\tau}-1)t}{\epsilon\tau}$.
\label{uLvLprop}
\end{my_prop}
\begin{proof}
By the definition of $u_L$ and $v_L$ given in (\ref{uL}) and (\ref{vL}) and the assumption that $u_0(x)=v_0(x)$ for $x\in[0,L]$, we can get their difference
\begin{equation*}
u_L(x,t)-v_L(x,t)=
c_1(t)\left(e^{-\frac{x}{\epsilon\sqrt{\tau}}}-\phi_1(x)\right)+\left(U_{c_2}(L,t)-h(t)+e^{-\frac{t}{\epsilon\tau}}h(0)\right)\phi_2(x)
\end{equation*}
Combining Lemmas \ref{lemma_old4}(\ref{lemmaphi1}), \ref{lemma_old4}(\ref{lemmaphi2}), inequality (\ref{Uc2}), and $h(t)\equiv 0$, we have
\begin{equation}
\norm{u_L(\cdot,t)-v_L(\cdot,t)}_\infty  \le  E_{1;\epsilon,\tau}(t)  e^{-\frac{\lambda
L}{\epsilon\sqrt{\tau}}}+ E_{2;\epsilon,\tau}(t)
e^{-\frac{\lambda(L-L_0)}{\epsilon\sqrt{\tau}}}\label{uLvL}
\end{equation}
where
\begin{eqnarray}
E_{1;\epsilon,\tau}(t) =
|c_1(\cdot)|_\infty+a_\tau e^\frac{b_\tau t}{\epsilon\tau}
\quad\mathrm{and}\quad
E_{2;\epsilon,\tau}(t) =
c_{\tau}\frac{t}{\epsilon\tau}e^\frac{(b_{\tau}-1)t}{\epsilon\tau}.
\label{Edef}
\end{eqnarray}
\end{proof}

\begin{my_prop}
If $u_0(x)$ satisfies (\ref{u0}), and $E_{1;\epsilon,\tau}(t), E_{2;\epsilon,\tau}(t)$ are as in proposition \ref{uLvLprop}, then 
\begin{equation*}
\norm{u_L(\cdot,t)-v_L(\cdot,t)}_{H_{L,\epsilon,\tau}^1}\le
\sqrt{5L}\left(E_{1;\epsilon,\tau}(t)e^{-\frac{\lambda
L}{\epsilon\sqrt{\tau}}}+E_{2;\epsilon,\tau}(t)e^{-\frac{\lambda(L-L_0)}{\epsilon\sqrt{\tau}}}\right)
.
\end{equation*}
\label{uLvLH1prop}
\end{my_prop}
\begin{proof}
Because of the definition of $u_L$ and $v_L$ given in (\ref{uL}) and (\ref{vL}), Lemma \ref{lemma_old4}(\ref{lemmaphi2prime}) and inequality  (\ref{Uc2}), we have that
\begin{equation}
\begin{split}
\norm{(u_L(\cdot,t)-v_L(\cdot,t))_x}_\infty \le&
\left|c_1(t)\right|e^{-\frac{L}{\epsilon\sqrt{\tau}}}\left|\phi_2'(x)\right|+\left|U_{c_2}(L,t)\right|\left|\phi_2'(x)\right|\\
\le&\frac{2}{\epsilon\sqrt{\tau}}\left( E_{1;\epsilon,\tau}(t)  e^{-\frac{\lambda
L}{\epsilon\sqrt{\tau}}}+ E_{2;\epsilon,\tau}(t)
e^{-\frac{\lambda(L-L_0)}{\epsilon\sqrt{\tau}}} \right) . 
\end{split}
\label{uLvLprime}
\end{equation}
Now, combining (\ref{uLvL}) and (\ref{uLvLprime}), we obtain that
\begin{equation}
\begin{split}
\norm{u_L(\cdot,t)-v_L(\cdot,t)}_{H_{L,\epsilon,\tau}^1}=&\sqrt{\int_0^L\left|u_L-v_L\right|^2+\left|\epsilon\sqrt{\tau}\left(u_L-v_L\right)_x\right|^2\,dx}\\
\le&
\sqrt{5L}\left(E_{1;\epsilon,\tau}(t)e^{-\frac{\lambda
L}{\epsilon\sqrt{\tau}}}+E_{2;\epsilon,\tau}(t)e^{-\frac{\lambda(L-L_0)}{\epsilon\sqrt{\tau}}}\right) .
\end{split}
\label{uLvLH1}
\end{equation}
\end{proof}

\begin{my_prop}
If $u_0(x)$ satisfies (\ref{u0}), then 
\begin{equation*}
\norm{W(\cdot,t)}_{H_{L,\epsilon,\tau}^1} \le
\gamma_{1;\epsilon,\tau}(t)e^{-\frac{\lambda
L}{\epsilon\sqrt{\tau}}}+\gamma_{2;\epsilon,\tau}(t)e^{-\frac{\lambda
(L-L_0)}{\epsilon\sqrt{\tau}}}
\end{equation*}
where the coefficients are given by
\begin{equation}
\begin{split}
\gamma_{1;\epsilon,\tau}(t)& = 
e^\frac{(M+1)^2t}{2M\epsilon\sqrt{\tau}}\left(\frac{(M+1)^2\sqrt{\tau}}{2M}+1\right)\sqrt{L}
\left(
\frac{t}{\epsilon\tau}|c_1(\cdot)|_\infty+\frac{a_\tau}{b_\tau}(e^\frac{b_\tau t}{\epsilon\tau}-1)
\right)\\
\gamma_{2;\epsilon,\tau}(t)
&= e^\frac{(M+1)^2t}{2M\epsilon\sqrt{\tau}}\left(\frac{(M+1)^2\sqrt{\tau}}{2M}+1\right)\sqrt{L}
c_\tau\\
&\qquad\cdot
\left(
\frac{t}{\epsilon\tau(b_\tau-1)}e^\frac{(b_\tau-1)t}{\epsilon\tau}-\frac{1}{(b_\tau-1)^2}(e^\frac{(b_\tau-1)t}{\epsilon\tau}-1)
\right).
\end{split}
\label{gammadef}
\end{equation}
 \label{Wthm}
\end{my_prop}
\begin{proof}
Multiplying the governing equation of $W$ (\ref{Weq}) by $2W$, integrating over $[0,L]$, and
using integration by parts, we get
\begin{eqnarray*}
&&\frac{d}{dt}\int_0^L W^2+(\epsilon\sqrt{\tau}W_x)^2\,dx\\
&=&-\epsilon\int_0^L 2W_x^2\, dx +\int_0^L
2W_x\left(f(v)-f(u)\right)\, dx
+\frac{2}{\epsilon\tau}\int_0^L W(v_L-u_L)\, dx .
\end{eqnarray*}
Therefore, using the norm we defined earlier in (\ref{Ynorm}), and 
$f'(u)\le\frac{(M+1)^2}{2M}:=C$, we have
\begin{eqnarray*}
&&\frac{d}{dt} \norm{W(\cdot,t)}_{H_{L,\epsilon,\tau}^1}^2\\
&\le & 2\int_0^L|W_x||f'(\eta)||v-u|\,dx
+\frac{2\sqrt{L}}{\epsilon\tau}\norm{v_L-u_L}_\infty\norm{W(\cdot,t)}_{H_{L,\epsilon,\tau}^1}\\
&\le&2C\int_0^L|W_x|\left(|W|+\norm{v_L-u_L}_\infty\right)\,dx+\frac{2\sqrt{L}}{\epsilon\tau}\norm{v_L-u_L}_\infty\norm{W(\cdot,t)}_{H_{L,\epsilon,\tau}^1}\\
&\le& 
\frac{2C}{\epsilon\sqrt{\tau}}\left(\norm{W(\cdot,t)}_{H_{L,\epsilon,\tau}^1}^2+\norm{v_L-u_L}_\infty\sqrt{L}\norm{ 
W(\cdot,t)}_{H_{L,\epsilon,\tau}^1}\right)\\
&&\qquad+\frac{2\sqrt{L}}{\epsilon\tau}\norm{v_L-u_L}_\infty\norm{W(\cdot,t)}_{H_{L,\epsilon,\tau}^1}\\
&=& \frac{2C}{\epsilon\sqrt{\tau}}\norm{ W(\cdot,t)}_{H_{L,\epsilon,\tau}^1}^2+\left( 
\frac{2C}{\epsilon\sqrt{\tau}}+\frac{2}{\epsilon\tau}\right)\sqrt{L}\norm{v_L-u_L}_\infty\norm{W(\cdot,t)}_{H_{L,\epsilon,\tau}^1} 
.
\end{eqnarray*}
Hence,
\begin{equation*}
\frac{d}{dt} \norm{ W(\cdot,t)}_{H_{L,\epsilon,\tau}^1}\le \frac{C}{\epsilon\sqrt{\tau}}\norm{W(\cdot,t)}
_{H_{L,\epsilon,\tau}^1}+\left(\frac{C}{\epsilon\sqrt{\tau}}+\frac{1}{\epsilon\tau}\right)\sqrt{L}\norm{v_L-u_L}_\infty 
.
\end{equation*}
By Gronwall's inequality and (\ref{uLvL})
\begin{eqnarray*}
&&\norm{ W(\cdot,t)}_{H_{L,\epsilon,\tau}^1}  \\
&\le&
\int_0^t\left(\frac{C}{\epsilon\sqrt{\tau}}+\frac{1}{\epsilon\tau}\right)\sqrt{L}\norm{v_L-u_L}_\infty 
e^\frac{C(t-s)}{\epsilon\sqrt{\tau}}\,ds\\
&\le& 
e^\frac{Ct}{\epsilon\sqrt{\tau}}\left(\frac{C}{\epsilon\sqrt{\tau}}+\frac{1}{\epsilon\tau}\right)\sqrt{L}\int_0^t
E_{1;\epsilon,\tau}(s)e^{-\frac{\lambda L}{\epsilon\sqrt{\tau}}} + E_{2;\epsilon,\tau}(s)e^{-\frac{\lambda(L-L_0)}{\epsilon\sqrt{\tau}}}ds\\
&\le&
~~\left(e^\frac{Ct}{\epsilon\sqrt{\tau}}\left(\frac{C}{\epsilon\sqrt{\tau}}+\frac{1}{\epsilon\tau}\right)\sqrt{L}\int_0^t
E_{1;\epsilon,\tau}(s)\,ds\right) e^{-\frac{\lambda L}{\epsilon\sqrt{\tau}}}\\
&&+
\left(e^\frac{Ct}{\epsilon\sqrt{\tau}}\left(\frac{C}{\epsilon\sqrt{\tau}}+\frac{1}{\epsilon\tau}\right)\sqrt{L}\int_0^tE_{2;\epsilon,\tau}(s)\,ds\right)
e^{-\frac{\lambda (L-L_0)}{\epsilon\sqrt{\tau}}} \\
&\le&~~ 
e^\frac{Ct}{\epsilon\sqrt{\tau}}\left(\frac{C}{\epsilon\sqrt{\tau}}+\frac{1}{\epsilon\tau}\right)\sqrt{L}
\left(
t|c_1(\cdot)|_\infty+\frac{a_\tau\epsilon\tau}{b_\tau}(e^\frac{b_\tau t}{\epsilon\tau}-1)
\right) e^{-\frac{\lambda L}{\epsilon\sqrt{\tau}}}\\
&&+
e^\frac{Ct}{\epsilon\sqrt{\tau}}\left(\frac{C}{\epsilon\sqrt{\tau}}+\frac{1}{\epsilon\tau}\right)\sqrt{L}
\frac{c_\tau}{\epsilon\tau}
\left(
\frac{\epsilon\tau}{b_\tau-1}te^\frac{(b_\tau-1)t}{\epsilon\tau}-(\frac{\epsilon\tau}{b_\tau-1})^2(e^\frac{(b_\tau-1)t}{\epsilon\tau}-1)
\right)e^{-\frac{\lambda (L-L_0)}{\epsilon\sqrt{\tau}}}.
\end{eqnarray*}
Hence
\begin{equation*}
\norm{ W(\cdot,t)}_{H_{L,\epsilon,\tau}^1}  \le
\gamma_{1;\epsilon,\tau}(t)e^{-\frac{\lambda
L}{\epsilon\sqrt{\tau}}} +
\gamma_{2;\epsilon,\tau}(t)e^{-\frac{\lambda
(L-L_0)}{\epsilon\sqrt{\tau}}}
\end{equation*}
where $\gamma_{1;\epsilon,\tau}(t)$ and $\gamma_{2;\epsilon,\tau}(t)$ are given in \eqref{gammadef}.
\end{proof}

Now we are in the position to 
prove the main theorem of this section.

\begin{my_thm}\label{main}
If $u_0(x)$ satisfies
\begin{eqnarray*}
u_0(x)=\left\{
\begin{array}{llll}
 C_u & & &x\in [0,L_0]\\
0& & &x > L_0
\end{array}
\right.
\end{eqnarray*}
where $L_0<L$ and $C_u,$  
%
are positive constants, 
and $E_{1;\epsilon,\tau}(t),
E_{2;\epsilon,\tau}(t),\gamma_{1;\epsilon,\tau}(t),
\gamma_{2;\epsilon,\tau}(t) $ are as in \eqref{Edef} and \eqref{gammadef}
, then 
\begin{equation}
\norm{ u(\cdot,t)-v(\cdot,t)}_{H_{L,\epsilon,\tau}^1}\le D_{1;\epsilon,\tau}(t)
e^{-\frac{\lambda L}{\epsilon\sqrt{\tau}}}+D_{2;\epsilon,\tau}(t)
e^{-\frac{\lambda (L-L_0)}{\epsilon\sqrt{\tau}}}
\label{final_est}
\end{equation}
for some $0<\lambda<1$, and
\begin{equation*}
D_{1;\epsilon,\tau}(t) =
\gamma_{1;\epsilon,\tau}(t)+\sqrt{5L}E_{1;\epsilon,\tau}(t),\qquad
D_{2;\epsilon,\tau}(t) =
\gamma_{2;\epsilon,\tau}(t)+\sqrt{5L}E_{2;\epsilon,\tau}(t).
\end{equation*}
\end{my_thm}
\begin{proof}[Proof of the Main Theorem]
\begin{equation*}
\begin{split}
\norm{ u(\cdot,t)-v(\cdot,t)}_{H_{L,\epsilon,\tau}^1} &\le
\norm{ W(\cdot,t) }_{H_{L,\epsilon,\tau}^1} + \norm{
v_L(\cdot,t)-u_L(\cdot,t)}_{H_{L,\epsilon,\tau}^1}\\
&=D_{1;\epsilon,\tau}(t)  e^{-\frac{\lambda
L}{\epsilon\sqrt{\tau}}}+D_{2;\epsilon,\tau}(t) e^{-\frac{\lambda
(L-L_0)}{\epsilon\sqrt{\tau}}}
\end{split}
\end{equation*}
where
\begin{equation*}
\begin{split}
D_{1;\epsilon,\tau}(t) =&
\gamma_{1;\epsilon,\tau}(t)+\sqrt{5L}E_{1;\epsilon,\tau}(t)\\
=&e^\frac{(M+1)^2t}{2M\epsilon\sqrt{\tau}}\left(\frac{(M+1)^2\sqrt{\tau}}{2M}+1\right)\sqrt{L}
\left(
\frac{t}{\epsilon\tau}|c_1(\cdot)|_\infty+\frac{a_\tau}{b_\tau}(e^\frac{b_\tau t}{\epsilon\tau}-1)
\right)\\
&+\sqrt{5L}(|c(\cdot)|_\infty+a_\tau e^\frac{b_\tau t}{\epsilon\tau}),\\
D_{2;\epsilon,\tau}(t) =& \gamma_{2;\epsilon,\tau}(t)+\sqrt{5L}E_{2;\epsilon,\tau}(t)\\
=&e^\frac{(M+1)^2t}{2M\epsilon\sqrt{\tau}}\left(\frac{(M+1)^2\sqrt{\tau}}{2M}+1\right)\sqrt{L}
c_\tau\cdot\\
&\hspace{3cm}\cdot
\left(
\frac{t}{\epsilon\tau(b_\tau-1)}e^\frac{(b_\tau-1)t}{\epsilon\tau}-\frac{1}{(b_\tau-1)^2}(e^\frac{(b_\tau-1)t}{\epsilon\tau}-1)
\right)\\
&+\sqrt{5L}c_\tau\frac{t}{\epsilon\tau} e^\frac{(b_\tau-1)t}{\epsilon\tau}.
\end{split}
\end{equation*}
\end{proof}

This result gives that $\norm{u(\cdot,t)-v(\cdot,t)}_{H^1_{L,\epsilon,\tau}}$ 
exponentially delays in $L$.
This theorem shows that if $\frac{\lambda L}{\epsilon\sqrt{\tau}}$
and $\frac{\lambda (L-L_0)}{\epsilon\sqrt{\tau}}$ converge to
infinity, then the solution $v(x,t)$ of the finite interval
problem converges to the solution $u(x,t)$ of the
half line problem in the sense of
$\norm{\cdot}_{H_{L,\epsilon,\tau}^1}$. This can be achieved
either by letting $L\rightarrow\infty$ or $\epsilon\rightarrow 0$.
For example, in the extreme case, $\epsilon=0$, the half line
problem (\ref{quarter_plane_0}) becomes hyperbolic and the domain of
dependence is finite, so, certainly, one only need to consider the
finite interval problem. 
This is consistent with
the main theorem in the sense that for a fixed final time $t$, if $\lambda L>b_\tau t$ and 
$\lambda(L-L_0)>(b_\tau-1)t$, i.e., $L>\max(\frac{b_\tau t}{\lambda}, \frac{(b_\tau -1)t}{\lambda})$, then
$\norm{ u(\cdot,t)-v(\cdot,t)
}_{H_{L,\epsilon,\tau}^1}\le D_{1;\epsilon,\tau}(t)
e^{-\frac{\lambda L}{\epsilon\sqrt{\tau}}}+D_{2;\epsilon,\tau}(t)
e^{-\frac{\lambda (L-L_0)}{\epsilon\sqrt{\tau}}}\rightarrow 0$ as
$\epsilon\rightarrow 0$. 
Theorem \ref{main} gives a theoretical
justification for using the solution of the finite interval problem to approximate the solution of the half
line problem with appropriate choice of $L$ and $\epsilon$. Hence
in the next chapter, the numerical scheme designed to solve the
MBL equation (\ref{MBL}) is given for finite interval problem.

\section{Numerical schemes}
\label{central}
%
To numerically solve the MBL equation \eqref{MBL}, 
We first collect all the terms with time derivative and  rewrite
MBL equation (\ref{MBL}) as
\begin{equation}
(u-\epsilon^2\tau u_{xx})_t +(f(u))_x = \epsilon
u_{xx}.\label{evolution}
\end{equation}
By letting
\begin{equation}
w=u-\epsilon^2\tau u_{xx} \quad\Longleftrightarrow\quad u=(I-\epsilon^2\tau\partial_{xx})^{-1}w,\label{w_def}
\end{equation}
MBL equation (\ref{evolution}) can be written as
\begin{eqnarray}
w_t +(f(u))_x = \epsilon u_{xx}.\label{w_evolution}
\end{eqnarray}
Now, the new form of MBL equation (\ref{w_evolution}) can be viewed as a PDE in terms of $w$, and the occurrence of $u$ can be  recovered by (\ref{w_def}).   Equation (\ref{w_evolution}) can be formally viewed as
\begin{eqnarray}
w_t+(f((I-\epsilon^2\tau\partial_{xx})^{-1}w))_x=\epsilon((I-\epsilon^2\tau\partial_{xx})^{-1}w)_{xx},
\end{eqnarray}
which is a balance law in term of $w$. 
We adopt numerical schemes originally designed
for hyperbolic equations to solve the MBL equation (\ref{evolution}),
which is of pseudo-parabolic type. The local discontinuous Galerkin method has
been applied to solve equations involving mixed derivatives $u_{xxt}$ term \cite{Shu1,Shu2}.
To the best knowledge of the authors, the central schemes have not been applied to
solve equations of this kind. The main advantage of the central schemes is the simplicity.
``the direction of the wind`` is not required to be identified, and hence
the field-by-field decomposition can be avoided. In this chapter, we demonstrate
how to apply the central schemes to solve the MBL equation (\ref{evolution}).
\subsection{Second-order schemes}
\label{2ndorderscheme} In this section, we show how to apply the classical second
order central schemes \cite{NT} originally designed for hyperbolic
conservation laws to numerically solve the MBL equation
(\ref{MBL}), which is of pseudo-parabolic type. 
To solve (\ref{w_evolution}), we modify the central scheme given in \cite{NT}.
As in \cite{NT}, at each time level, we first reconstruct a
piecewise linear approximation of the form
\begin{eqnarray}
L_j(x,t)=w_j(t)+(x-x_j)\frac{w_j'}{\Delta x}\;, \qquad
x_{j-\frac{1}{2}}\le x \le x_{j+\frac{1}{2}}. \label{linear}
\end{eqnarray}
Second-order accuracy is guaranteed if the so-called vector of
numerical derivative $\frac{w_j'}{\Delta x}$, which will be given later, satisfies
\begin{eqnarray}
\frac{w_j'}{\Delta x}=\frac{\partial w(x_j,t)}{\partial
x}+O(\Delta x). \label{wprime}
\end{eqnarray}
We denote the staggered piecewise-constant functions
$\bar{w}_{j+\frac{1}{2}}(t)$ as
\begin{eqnarray}
\bar{w}_{j+\frac{1}{2}}(t) = \frac{1}{\Delta
x}\int_{x_j}^{x_{j+1}}w(x,t)\, dx .\label{staggered_grid}
\end{eqnarray}
Evolve the piecewise linear interplant (\ref{linear}) by
integrating (\ref{w_evolution}) over
$[x_j,x_{j+1}]\times[t,t+\Delta t]$
\begin{equation}
\begin{split}
\bar{w}_{j+\frac{1}{2}}(t+\Delta t) =&
\bar{w}_{j+\frac{1}{2}}(t)
 \\
& -\frac{1}{\Delta x}\left[ \int_{t}^{t+\Delta
t}f(u(x_{j+1},s))\,ds - \int_{t}^{t+\Delta
t}f(u(x_{j},s))\,ds \right]  \\
& +\frac{\epsilon}{\Delta x} \left[\int_t^{t+\Delta
t}\int_{x_j}^{x_{j+1}} \frac{\partial^2 u(x,s)}{\partial x^2}
\,dx\,ds\right]. 
\end{split}\label{evolution_int}
\end{equation}
We calculate each term on the right hand side of \eqref{evolution_int} below.
For $\bar{w}_{j+\frac{1}{2}}(t)$, applying the definition of $L_j(x,t)$ and $L_{j+1}(x,t)$ given 
in (\ref{linear}) to \eqref{staggered_grid}, we have that
\begin{equation}
\begin{split}
\bar{w}_{j+\frac{1}{2}}(t)&=\frac{1}{\Delta
x}\int_{x_j}^{x_{j+\frac{1}{2}}}L_j(x,t) \, dx +\frac{1}{\Delta
x}\int_{x_{j+\frac{1}{2}}}^{x_{j+1}}L_{j+1}(x,t) \, dx \\
&= \frac{1}{2}(w_j(t)+w_{j+1}(t))+\frac{1}{8}(w'_j-w'_{j+1}).
\end{split}
\label{w_ave}
\end{equation}
The middle two integrands can be approximated by the
midpoint rule
\begin{equation}
\begin{split}
\int_{t}^{t+\Delta t}f(u(x_j,s))\, ds &= f(u(x_j,t+\frac{\Delta
t}{2}))\Delta t + O(\Delta t^3)
\\
\int_{t}^{t+\Delta t}f(u(x_{j+1},s))\, ds &=
f(u(x_{j+1},t+\frac{\Delta t}{2}))\Delta t + O(\Delta t^3)
\end{split}
\label{integratef}
\end{equation}
if the CFL condition
\begin{equation*}
\lambda\cdot \max_{x_j\le x\le
x_{j+1}}\left|\frac{\partial f(u(w(x,t)))}{\partial w}\right|<\frac{1}{2},\qquad \mathrm{where}\quad
\lambda=\frac{\Delta t}{\Delta x}
\end{equation*}
is met. 
For MBL equation \eqref{w_evolution}, we have that at $t>0$,
\begin{equation*}
 u-\epsilon^2\tau u_{xx}=w,\qquad u(0)=w(0),\qquad u(L)=w(L).
\end{equation*}
Let $v(x)=\frac{(L-x)w(0)+x w(L)}{L}$, then
\begin{equation*}
 u(x)=[Iw](x)=v(x)+\frac{1}{L}\int_0^L[w(y)-v(y)]\,K(x,y)\,dy
\end{equation*}
where
\begin{equation*}
 K(x,y)=\sum_{k=1}^{\infty}\frac{\sin(\frac{k\pi x}{L})\sin(\frac{k\pi y}{L})}{1+(\frac{k\pi}{L})^2\epsilon^2\tau}.
\end{equation*}
Hence the eigenvalues for $I$ are
\begin{equation*}
 \lambda_k=\frac{1}{1+(\frac{k\pi}{L})^2\epsilon^2\tau}\le 1, \qquad k=1,2,3\dots
\end{equation*}
Therefore, the CFL condition is
\begin{equation*}
\frac{\Delta t}{\Delta x}\cdot\max_{x_j\le x\le x_{j+1}}\left|\frac{\partial f(u(w(x,t)))}{\partial w}\right| 
=
\frac{\Delta t}{\Delta x}\cdot\max_{\substack{x_j\le x\le x_{j+1}\\k=1,2,3\dots}}\left|\frac{\partial f(u(x,t))}{\partial u}\right|\cdot\lambda_k
\le
\frac{\Delta t}{\Delta x}\cdot 2.2
<\frac{1}{2}
\end{equation*}
In the numerical computations in chapter \ref{results}, we chose $\frac{\Delta t}{\Delta x}=0.1$.
%
%
%
%
In \eqref{integratef}, to estimate $u(\cdot,t+\frac{\Delta t}{2})$'s,
we use Taylor expansion and the conservation law
(\ref{w_evolution}):
\begin{equation}
\begin{split}
w(x_j,t+\frac{\Delta t}{2}) &=
w_j(t)+\frac{\partial w}{\partial t}\frac{\Delta t}{2} + \mathcal{O}(\Delta t^2) \\
&=w_j(t)+(\epsilon\frac{\partial^2 u}{\partial x^2}-\frac{\partial f}{\partial x})\frac{\Delta t}{2} + \mathcal{O}(\Delta t^2) \\
&=w_j(t)+(\epsilon\Delta x\,D^2\,u_j-f_j')\frac{\lambda}{2}, 
\end{split}
\label{halftime}
\end{equation}
where $D$ is the discrete central difference operator
\begin{equation*}
D^2 u_j = \frac{u_{j-1} - 2u_j + u_{j+1}}{\Delta x^2},
\end{equation*}
and the second-order accuracy is met if
\begin{eqnarray}
\frac{f_j'}{\Delta x}=\frac{\partial f(u(x_j,t))}{\partial
x}+\mathcal{O}(\Delta x). \label{fprime}
\end{eqnarray}
%
%
The choices for \{$w_j'$\} in \eqref{wprime} and
\{$f_j'$\} in \eqref{fprime} can be found in \cite{NT}, 
and we chose
\begin{eqnarray}
w_j'=MM\{\Delta w_{j+\frac{1}{2}}, \Delta
w_{j-\frac{1}{2}}\}\;,\qquad f_j'=MM\{\Delta f_{j+\frac{1}{2}},
\Delta f_{j-\frac{1}{2}}\}\label{standardvprime}
\end{eqnarray}
where
$MM\{x,y\}=\mathrm{minmod}(x,y)=\frac{1}{2}(\mathrm{sgn}(x)+\mathrm{sgn}(y))\cdot
\mathrm{Min}(|x|,|y|)$ and  $\Delta w_{j+\frac{1}{2}}=w_{j+1}-w_j$.
Combining (\ref{evolution_int})-(\ref{integratef}), we obtain
\begin{equation}
\begin{split}
\bar{w}_{j+\frac{1}{2}}(t+\Delta t) =&
\bar{w}_{j+\frac{1}{2}}(t)\\
&-\lambda[f(u_{j+1}(t+\frac{\Delta t}{2})-f(u_j(t+\frac{\Delta
t}{2}))] \\
& +\frac{\epsilon}{\Delta x} \left[\int_t^{t+\Delta
t}\int_{x_j}^{x_{j+1}} \frac{\partial^2 u(x,s)}{\partial x^2}
\,dx\,ds\right]
.
\end{split}
\label{staggered11}
\end{equation}
%
%
%
%
%
Next, we will re-write \eqref{staggered11} in terms of $u$.
$(\overline{u_{xx}})_{j+\frac{1}{2}}$ is approximated as 
\begin{equation*}
(\overline{u_{xx}})_{j+\frac{1}{2}}=\frac{1}{\Delta x}\int^{x_{j+1}}_{x_j} u_{xx}\,dx
=\frac{1}{\Delta x}(u_x(x_{j+1},t)-u_x(x_j,t)),
\end{equation*}
and using the cell averages, it becomes
\begin{equation}
\begin{split}
(\overline{u_{xx}})_{j+\frac{1}{2}}&=\frac{1}{\Delta x}\left( \frac{\bar{u}_{j+3/2} - \bar{u}_{j+1/2}}{\Delta x} - 
\frac{\bar{u}_{j+1/2} - \bar{u}_{j-1/2}}{\Delta x}\right)\\
&=\frac{\bar{u}_{j+3/2} - 2\bar{u}_{j+1/2}
+ \bar{u}_{j-1/2}}{(\Delta x)^2}\\
&=D^2\bar{u}_{j+\frac{1}{2}}.
\end{split}
\label{uxxbar}
\end{equation}
Notice that the linear interpolation  (similar to \eqref{linear})
\begin{equation*}
\tilde{L}_{j+\frac{1}{2}}(x,t+\Delta t) = u_{j+\frac{1}{2}}(t+\Delta t)+(x-x_{j+\frac{1}{2}})\frac{u'_{j+\frac{1}{2}}}{\Delta x}
\quad  \text{for} \quad  x_j\le x\le x_{j+1}
\end{equation*}
 and the cell average definition (similar to \eqref{staggered_grid}) 
\begin{equation*}
\bar{u}_{j+\frac{1}{2}}(t+\Delta t)=\frac{1}{\Delta t}\int_{x_j}^{x_{j+1}}u(x,t+\Delta t)\,dx
\end{equation*}
ensure that 
\begin{equation*} 
\bar{u}_{j+\frac{1}{2}}(t+\Delta t) = u_{j+\frac{1}{2}}(t+\Delta t),
\end{equation*}
and the convertion between $u$ and $w$ is done using the following relation
\begin{equation}
 (I-\epsilon^2\tau\,D^2)u=w. \label{convertion}
\end{equation}
Hence re-writting \eqref{staggered11} in terms of $u$ gives the staggered central scheme 
\begin{equation}
\begin{split}
(I-\epsilon^2\tau\,D^2)u_{j+\frac{1}{2}}(t+\Delta t)
&=(I-\epsilon^2\tau\,D^2)\bar{u}_{j+\frac{1}{2}}(t)\\
&-\lambda[f(u_{j+1}(t+\frac{\Delta t}{2})-f(u_j(t+\frac{\Delta
t}{2}))] \\
& +\frac{\epsilon}{\Delta x} \left[\int_t^{t+\Delta
t}\int_{x_j}^{x_{j+1}} \frac{\partial^2 u(x,s)}{\partial x^2}
\,dx\,ds\right].
\end{split}
\label{staggeredu}
\end{equation}
%
%
%
%
%
%
We will focus on the last integral in
(\ref{staggeredu}). There are many ways to numerically calculate
this integral. 
We will show two ways to do this in the following
two subsections, both of them achieve second order accuracy.
\subsubsection{Trapezoid Scheme}
\label{trapzoidsec} In this scheme, we use the notion
(\ref{staggered_grid}) and the
trapezoid rule to calculate the integral numerically as follows:
\begin{equation}
\begin{split}
\int_{t}^{t+\Delta t}\int_{x_j}^{x_{j+1}}\frac{\partial^2
u(x,s)}{\partial x^2}\,dx\,ds&=\Delta x\int_{t}^{t+\Delta
t}(\overline{u_{xx}})_{j+\frac{1}{2}}(s)\,ds\\
&=\frac{\Delta x\Delta
t}{2}\left((\overline{u_{xx}})_{j+\frac{1}{2}}(t)+(\overline{u_{xx}})_{j+\frac{1}{2}}(t+\Delta
t))\right)
\end{split}
\end{equation}
with $\mathcal{O}(\Delta t^3)$ error. Combining with
\eqref{uxxbar} and (\ref{staggeredu}), we can get the trapezoid scheme
%
%
%
\begin{equation}
\begin{split}
\left(I-(\epsilon^2\tau+\frac{\epsilon\Delta
t}{2})D^2\right)u_{j+\frac{1}{2}}(t+\Delta t) =\left(I-(
\epsilon^2\tau-\frac{\epsilon\Delta t}{2})D^2\right)
\bar{u}_{j+\frac{1}{2}}(t)
\\
-\lambda\left[f(u_{j+1}(t+\frac{\Delta
t}{2}))-f(u_{j}(t+\frac{\Delta t}{2}))\right]
\end{split}\label{trapzoid}
\end{equation}
The flow chart of the trapezoid scheme is given in \eqref{Flow:trapezoid}
\begin{equation}
\xymatrix@R=0.1pt{
&&\bar{w}_{j+\frac{1}{2}}(t)\ar[r]^{\eqref{convertion}}&\bar{u}_{j+\frac{1}{2}}(t)\ar[rd]^{\eqref{trapzoid}}&\\
u_j(t)\ar[r]^{\eqref{convertion}} &  w_j(t)\ar[ru]^{\eqref{w_ave}}\ar[rd]_{\eqref{halftime}}&&&u_{j+\frac{1}{2}}(t+\Delta t)\\
&&w_j(t+\frac{\Delta t}{2})\ar[r]^{\eqref{convertion}}&u_j(t+\frac{\Delta t}{2})\ar[ru]_{\eqref{trapzoid}}
}
\label{Flow:trapezoid}
\end{equation}


%
%
%
%

\subsubsection{Midpoint Scheme}
 In this scheme, we use the notion
(\ref{staggered_grid}) and the
midpoint rule to calculate the integral numerically as follows:
\begin{equation*}
\begin{split}
\int_{t}^{t+\Delta t}\int_{x_j}^{x_{j+1}}\frac{\partial^2
u(x,s)}{\partial x^2}\,dx\,ds&=\Delta x\int_{t}^{t+\Delta
t}(\overline{u_{xx}})_{j+\frac{1}{2}}(s)\,ds\nonumber\\
&=\Delta x\Delta t(\overline{u_{xx}})_{j+\frac{1}{2}}(t+\frac{\Delta t}{2})
\end{split}
\end{equation*}
Combining with
\eqref{uxxbar} and (\ref{staggeredu}), we can get the midpoint scheme
%
%
%
\begin{equation}
\begin{split}
(I-\epsilon^2\tau\, D^2){u}_{j+\frac{1}{2}}(t+\Delta t)  =&
\bar{w}_{j+\frac{1}{2}}(t)\\
&-\lambda[f(u_{j+1}(t+\frac{\Delta t}{2})-f(u_j(t+\frac{\Delta
t}{2}))]  \\
& +\epsilon\Delta t
D^2\bar{u}_{j+\frac{1}{2}}(t+\frac{\Delta t}{2})
\end{split}
\label{midpt}
\end{equation}
The flow chart of the midpoint scheme is given in \eqref{Flow:midpoint}
\begin{equation}
\xymatrix{
&&\bar{w}_{j+\frac{1}{2}}(t)\ar[drr]^{\eqref{midpt}}&&\\
u_j(t)\ar[r]^{\eqref{convertion}} & w_j(t)\ar[ur]^{\eqref{w_ave}}\ar[dr]_{\eqref{halftime}}&
\bar{w}_{j+\frac{1}{2}}(t+\frac{\Delta t}{2})\ar[r]^{\eqref{convertion}}&\bar{u}_{j+\frac{1}{2}}(t+\frac{\Delta t}{2})\ar[r]_{\eqref{midpt}}&u_{j+\frac{1}{2}}(t+\Delta t)\\
&&w_j(t+\frac{\Delta t}{2})\ar[u]_{\eqref{w_ave}}\ar[r]^{\eqref{convertion}}&
u_j(t+\frac{\Delta t}{2})\ar[ur]_{\eqref{midpt}}&&
}
\label{Flow:midpoint}
\end{equation}

\subsection{A third order semi-discrete  scheme}
\label{3rdorderscheme} 
Similarly, we can extend the third order scheme to solve MBL equation (\ref{MBL}), however, it is more involved. But the third order
semi-discrete central scheme proposed in \cite{KL} can be extended to solve the
MBL equation in a straightforward manner. In order to make the paper self-contained, we include the formulation below.
\begin{eqnarray*}
\frac{d\bar{w}_j}{dt} = -\frac{H_{j+1/2}(t) - H_{j-1/2}(t)}{\Delta x} + \epsilon Q_j(t) \\
\end{eqnarray*}
where $\bar{w}(x,t)$ denotes the cell average of $w$
\begin{eqnarray*}
\bar{w}_j(t)=\frac{1}{\Delta x}\int_{x_{j-1/2}}^{x_{j+1/2}}w(x,t)\,dx,
\end{eqnarray*}
 $H_{j+1/2}(t)$ is the numerical convection flux and  $Q_j(t)$ is a high-order approximation to the diffusion term $u_{xx}$.  
\begin{eqnarray*}
H_{j+1/2}(t) = \frac{f(u^+_{j+1/2}(t)) +
f(u^-_{j+1/2}(t))}{2}-\frac{a_{j+1/2}(t)}{2}\left[w^+_{j+1/2}(t)-w^-_{j+1/2}(t)\right]
\end{eqnarray*}
where $u^-_{j+1/2}(t),u^+_{j+1/2}(t) $ denote the left and right intermediate values of $u(x,t^n)$ at $x_{j+1/2}$, and their values are converted from the $w^-_{j+1/2}(t),w^+_{j+1/2}(t) $ using (\ref{w_def}). The way to calculate  $w^-_{j+1/2}(t)$, $w^+_{j+1/2}(t) $ and $a_{j+1/2}(t)$ is
\begin{equation*}
 \begin{split}
w^+_{j+1/2}(t) &=  A_{j+1} - \frac{\Delta x}{2}B_{j+1}
+\frac{(\Delta x)^2}{8}C_{j+1},\\
w^-_{j+1/2}(t) &=  A_{j} + \frac{\Delta x}{2}B_{j}
+\frac{(\Delta x)^2}{8}C_{j},\\
a_{j+1/2}(t) &= \mathrm{max}\left\{\frac{\partial f}{\partial u}(u_{j+1/2}^-(t)),\frac{\partial
f}{\partial u}(u_{j+1/2}^+(t))\right\} ,
\end{split}
\end{equation*}
where
\begin{equation*}
 \begin{split}
A_j &= \bar{w}_j^n-\frac{w_C}{12}(\bar{w}_{j+1}^n-2\bar{w}_j^n+\bar{w}_{j-1}^n), \\
B_j &=\frac{1}{\Delta x}\left[w_R(\bar{w}_{j+1}^n-\bar{w}_j^n)+w_C\frac{\bar{w}_{j+1}^n-\bar{w}_{j-1}^n}{2}+w_L(\bar{w}_j^n-\bar{w}_{j-1}^n)\right],\\
C_j &=2w_C\frac{\bar{w}_{j-1}^n-2\bar{w}_j^n+\bar{w}_{j+1}^n}{\Delta x^2},\\
w_i &= \frac{\alpha_i}{\sum_m \alpha_m} \qquad
\alpha_i=\frac{c_i}{(\epsilon_0+IS_i)^p}, \qquad i,m\in\{C,R,L\}\\
c_L&=c_R=1/4, \qquad c_C=1/2,\qquad \epsilon_0=10^{-6}, \qquad
p=2,\\
IS_L &= (\bar{w}_j^n-\bar{w}_{j-1}^n)^2 ,\qquad IS_R = (\bar{w}_{j+1}^n-\bar{w}_j^n)^2,\\
IS_C&=\frac{13}{3}(\bar{w}_{j+1}^n-2\bar{w}_j^n+\bar{w}_{j-1}^n)^2+\frac{1}{4}(\bar{w}_{j+1}^n-\bar{w}_{j-1}^n)^2.\\
\end{split}
\end{equation*}
The diffusion $u_{xx}$ is approximated using the following fourth-order central differencing form
\begin{eqnarray}
 Q_j(t) =\frac{-u_{j-2}+16u_{j-1}-30u_j+16u_{j+1}-u_{j+2}}
{12\Delta x^2}.
\label{Qj}
\end{eqnarray}
%
%
The unique feature of this scheme is that the discretization is done in space 
first, and then  the time  evolution equation can be solved as a system of 
ordinary differential equations using any ODE solver of third order or higher. 
In this paper, we simply use the standard fourth order Runge-Kutta methods.
Notice that to achieve the third order accuracy, the linear solver that 
converts $u$ from $w$ using (\ref{w_def}) need also to be high order,
and \eqref{Qj} is used to discretize $u_{xx}$ in our convertion.


\section{Computational results}
\label{results}
In this section, we show the numerical solutions to the MBL equation
\begin{eqnarray}
u_t+(f(u))_x=\epsilon u_{xx}+\epsilon^2\tau u_{xxt}\label{MBL_new}
\end{eqnarray}
with the initial condition
\begin{eqnarray}
u_0(x) = \left\{
\begin{array}{lll}
u_B & \qquad\mathrm{if} & x= 0\\
0 & \qquad\mathrm{if} & x>0
\end{array}
\right.
\label{halfline}
\end{eqnarray}
and the Dirichlet boundary condition.
%

\begin{figure}[htbp]
\subfiguretopcaptrue
\begin{center}
\subfigure[$t = 0.1$]{\label{diff_L_t=0.1}\includegraphics[width=0.325\textwidth]{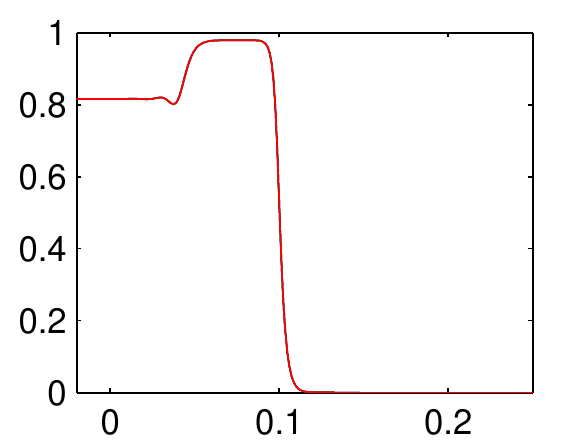}}
\subfigure[$t = 0.5$]{\label{diff_L_t=0.5}\includegraphics[width=0.325\textwidth]{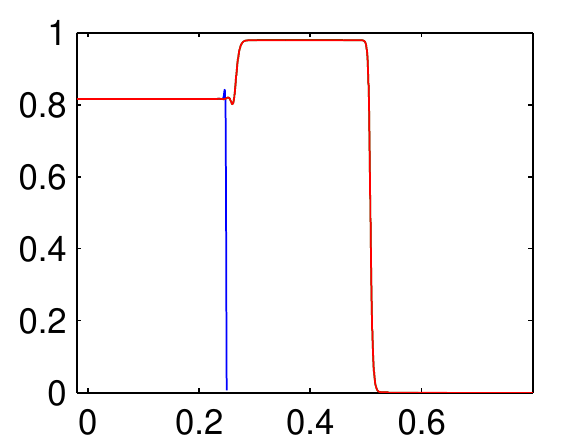}}
\subfigure[$t = 1$]{\label{diff_L_t=1}\includegraphics[width=0.325\textwidth]{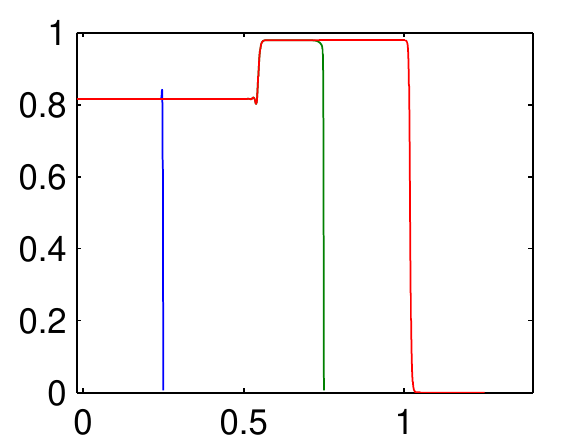}}
\caption{Numerical solutions of MBL \eqref{MBL_new}  at \subref{diff_L_t=0.1} $t = 0.1$, \subref{diff_L_t=0.5} $t = 0.5$,
\subref{diff_L_t=1} $t = 1$  using the trapezoid scheme \eqref{Flow:trapezoid}. `\textcolor{blue}{---}',  
`\textcolor{green}{---}',  `\textcolor{red}{---}' denote the numerical
solutions corresponding to computational domain $[0,L]$ with $L = 0.25$, $L = 0.75$ and $L = 1.25$ respectively.
The parameter values are  $\tau = 5$, $u_B = \alpha =\sqrt{\frac{2}{3}}$, $\epsilon = 0.001$, $\Delta x = \frac{\epsilon}{10}$,
$\Delta t = 0.1\Delta x$. }
\label{diff_L}
\end{center}
\end{figure}

Numerically, it is not practical to solve the half line problem \eqref{halfline}, and one
has to choose an appropriate computational domain.
Theorem \ref{main} in Chapter \ref{comparison} provides a theoretical bound for
the difference between the solution to the half line problem and that to
the finite interval problem.
However, the estimate \eqref{final_est} in Theorem \ref{main} includes time-dependent 
parameters $D_{1;\epsilon,\tau}(t)$ and $D_{2;\epsilon,\tau}(t)$, which cannot
be obtained analyticaly. 
Therefore, we numerically demonstrate how the computational domain size
affects the solution.
We choose $\tau = 5$, $u_B = \alpha =\sqrt{\frac{2}{3}}$ and
$\epsilon = 0.001$ as an example here.
Figure \ref{diff_L} shows the snapshot of the solutions at $t = 0.1$, $t = 0.5$ and $t = 1$ for
computational domain $[0,L]$ with $L = 0.25$, $L = 0.75$ and $L = 1.25$.
%
%
%

In Figure \ref{diff_L}\subref{diff_L_t=0.1}, $t = 0.1$, the leading shock is located at 
$\frac{f(\bar{u}_{\tau = 5})}{\bar{u}_{\tau = 5}} \times 0.1 = 1.02 \times 0.1 = 0.102$,
and $L = 0.25$, $L = 0.75$, $L = 1.25$ all exceed the leading shock location.
Hence all the three computational domains deliver visually indistinguishable results.
Whereas, in Figure \ref{diff_L}\subref{diff_L_t=0.5}, $t = 0.5$, the leading shock is
located at  $1.02 \times 0.5 = 0.51$, $L = 0.25$ is shorter than the computational domain
needed to capture this shock, hence the numerical solution halts at $x = 0.25$.
On the contrast, $L = 0.75$ and $L = 1.25$ are both large enough to capture this shock front.
Similarly, in Figure \ref{diff_L}\subref{diff_L_t=1}, $t = 1$, the leading shock is located at 
1.02. $L = 0.25 < 1.02$ and $L = 0.75 < 1.02$ both result in wrong solution profiles.
More specifically, both solutions halt at the boundary of the insufficient computational domain.
But $L = 1.25 > 1.02$ is large enough to capture the correct solution profile.

In the rest of this chapter, all the computational domains $[0,L]$ are therefore chosen based on the principle: 
$$L > \text{leading shock speed} \times \text{computational time}.$$
In addition, numerical solutions for larger $L$'s, for example, $L = 1.75$, $L = 2.5$, $L = 5$, $L = 10$ are also 
sought. For all these larger $L$'s, the numerical solutions are all consistent with that corresponding to $L = 1.25$ up to $t = 1$. This confirms 
that it is not necessary to take $L$ too much larger than 
{\tt leading shock speed} $\times$ {\tt computational time}.

%
To validate the order analysis given in chapter \ref{central} for various schemes proposed,
we first test the order of our schemes numerically
with a smooth initial
condition
\begin{eqnarray*}
u_0(x)=u_BH(x-5,5),
\end{eqnarray*}
where
\begin{eqnarray*}
H(x,\xi)=\left\{
\begin{array}{lll}
1 & \mathrm{if} & x<-\xi\\
1-\frac{1}{2}(1+\frac{x}{\xi}+\frac{1}{\pi}\sin(\frac{\pi x}{\xi})) & \mathrm{if} & -\xi\le x\le\xi\\
0&\mathrm{if} & x>\xi
\end{array}
\right. .
\end{eqnarray*}
The final time $T=1$ was employed, so that there was no shock created.
 $\epsilon$ in the MBL equation (\ref{MBL_new}) is taken to be 1,  $M$ is taken to be 2, and the computational 
interval is $[-10,20]$.
The $L_1, L_2, L_{\infty}$ order tests of the trapezoid scheme and the third order semi-discrete scheme with 
different parameter $\tau$ value and the initial condition $u_B$ are given in Tables \ref{trapezoid_2nd_order}, 
\ref{3rd_order}.
Table \ref{trapezoid_2nd_order} shows that the trapezoid rule achieved second order accuracy for all the tested 
cases in $L_1, L_2, L_{\infty}$ sense. Table \ref{3rd_order} shows that the semi-discrete scheme has the order of 
accuracy greater than 2.5 for all the cases, and exceeds 3 for some cases. This confirms the accuracy study given 
in sections \ref{trapzoidsec} and \ref{3rdorderscheme} respectively.

\begin{table}
\begin{center}
\begin{footnotesize}
\begin{tabular}{|l|l|l@{}l|l@{}l|l@{}l|}\hline
&&&&&&&\\
&N &$\norm{u_{\Delta x}-u_{\frac{\Delta x}{2}}}_1$ & order &$\norm{u_{\Delta x}-u_{\frac{\Delta x}{2}}}_2$ & 
order &$\norm{u_{\Delta x}-u_{\frac{\Delta x}{2}}}_\infty$ & order\\
&&&&&&&\\ \hline
   & 60&7.5416e-03  & -             &  2.5388e-03 &  -                  & 1.5960e-03&       -\\
$u_B=0.9$ & 120 &1.9684e-03 & 1.9379 &    6.7288e-04 &   1.9157  &   4.4066e-04  &   1.8568\\
$\tau=0.2$ & 240& 4.9891e-04 & 1.9802  &  1.7645e-04&   1.9311 &     1.2529e-04  &1.8144\\
 & 480& 1.2589e-04  & 1.9865&    4.5366e-05&   1.9596  &  3.3205e-05 &   1.9158\\\hline
    & 60& 8.0141e-03 & -&    2.6069e-03 & -&   1.4989e-03& -\\
 $u_B=0.9$ & 120 & 2.1502e-03  &  1.8981 &   7.0452e-04 &   1.8876&     4.2221e-04&   1.8279\\
 $\tau=1$ & 240&  5.5697e-04&    1.9488 &   1.8259e-04 &   1.9480  &  1.1283e-04 &   1.9038\\
  & 480& 1.4104e-04  &  1.9815 &    4.6109e-05 &  1.9855 &    2.8719e-05 &  1.9740\\\hline
   & 60&   1.3102e-02  &-&   4.1784e-03&-&   2.2411e-03 &-\\
$u_B=0.9$ & 120 & 3.6201e-03  & 1.8557  &   1.0994e-03  & 1.9263   & 6.1060e-04 &  1.8759\\
 $\tau=5$ & 240&  9.6737e-04  &  1.9039 &   2.8089e-04 &  1.9686 &   1.5667e-04  & 1.9625\\
  & 480&  2.5825e-04  &  1.9053 &   7.1250e-05&  1.9790 &    3.9286e-05 &  1.9956\\\hline
   & 60&     6.4427e-03  &- &    2.1578e-03 &- &    1.1682e-03  &-\\ 
$u_B=\alpha$ & 120 &    1.6611e-03  &   1.9555 &  5.7775e-04 &  1.9011  &   3.6447e-04  & 1.6804\\
 $\tau=0.2$ & 240&    4.3643e-04   &  1.9283  &   1.5215e-04  &   1.9250   & 1.0389e-04   & 1.8107\\
  & 480&   1.1223e-04  & 1.9593   &  3.9170e-05    & 1.9577  &  2.7629e-05  & 1.9109\\\hline
    & 60&  7.5867e-03  &-& 2.4101e-03  &-& 1.3364e-03&- \\
$u_B=\alpha$ & 120 &     2.0069e-03    &1.9185  &  6.4998e-04&  1.8906   &  3.7650e-04 &1.8277\\
 $\tau=1$ & 240&    5.1832e-04  & 1.9531 & 1.6801e-04&  1.9519   &  1.0062e-04 & 1.9037\\
  & 480&   1.3136e-04   &1.9803 &   4.2497e-05 &   1.9831   &  2.5599e-05  & 1.9748\\\hline
    & 60&    1.1959e-02 &-&   3.8026e-03 &-&  1.9938e-03  &-\\
$u_B=\alpha$ & 120 &    3.2940e-03  & 1.8602&  9.9527e-04 &  1.9338    & 5.4231e-04  &    1.8783\\
 $\tau=5$ & 240&   8.7736e-04 &  1.9086&   2.5358e-04 &   1.9727 &   1.3933e-04  &1.9606\\
  & 480&    2.3271e-04   & 1.9146&   6.4252e-05  &1.9806&     3.4967e-05  &1.9944\\\hline
    & 60&   5.7714e-03  &-& 1.9358e-03 &-&  1.0481e-03  &-\\
 $u_B=0.75$ & 120 &   1.5035e-03  &     1.9406 &  5.1617e-04   &1.9070&  2.8061e-04 &  1.9011\\
 $\tau=0.2$ & 240& 3.9299e-04  & 1.9357 &  1.3616e-04   &1.9225 & 7.9134e-05  &   1.8262\\
  & 480&   1.0063e-04   &  1.9655&  3.5080e-05  & 1.9566  & 2.1035e-05  &   1.9115\\\hline
   & 60&   7.1823e-03  &-&  2.2843e-03&-&   1.2069e-03  &-\\
  $u_B=0.75$ & 120 &  1.8963e-03  &  1.9213 &   6.1315e-04 &1.8974 &   3.4013e-03   & 1.8272\\
   $\tau=1$ & 240& 4.8284e-04    & 1.9736  & 1.5796e-04&    1.9567  & 9.0912e-04  &  1.9035\\
 & 480&     1.2093e-04  &  1.9974 &   3.9783e-05 &  1.9894 &  2.3121e-05  & 1.9753\\\hline
   & 60&    1.1042e-02   &-& 3.5020e-03  &-&  1.8299e-03  &-\\
  $u_B=0.75$ & 120 &   3.0287e-03  &1.8662  &  9.1181e-04& 1.9414&  4.8976e-04   & 1.9016\\
    $\tau=5$ & 240& 8.0111e-04  &  1.9186&  2.3118e-04  &1.9797  & 1.2593e-04 & 1.9595\\
 & 480&  2.1076e-04  &1.9264 &  5.8358e-05  &1.9860 &  3.1627e-05  & 1.9934\\\hline
 \end{tabular}
\end{footnotesize}
\end{center}
\caption{The accuracy test for the trapezoid scheme for the MBL equation (\ref{MBL_new}) with $\epsilon=1$ and 
$M=2$.}
\label{trapezoid_2nd_order}
\end{table}

\begin{table}
\begin{center}
\begin{footnotesize}
\begin{tabular}{|l|l|l@{}l|l@{}l|l@{}l|}\hline
&&&&&&&\\
&$N$ &$\norm{u_{\Delta x}-u_{\frac{\Delta x}{2}}}_1$ & order &$\norm{u_{\Delta x}-u_{\frac{\Delta 
x}{2}}}_2$ & order &$\norm{u_{\Delta x}-u_{\frac{\Delta x}{2}}}_\infty$ & order\\
&&&&&&&\\ \hline
   & 120&  2.6992e-03 &-&    1.1300e-03 &-&  7.2363e-04 &-\\
$u_B=0.9$ & 240 & 4.0403e-04 &  2.7400 &  1.7079e-04  &  2.7260  & 1.1283e-04  & 2.6811\\
$\tau=0.2$ & 480& 5.7504e-05&   2.8127 &  2.4624e-05  &  2.7941&   1.6242e-05 & 2.7963\\
 & 960&  8.4934e-06   &2.7592   &  3.0892e-06   &2.9948 &   1.7607e-06  &  3.2055\\\hline
   & 120&    4.7731e-03 & -&  2.0192e-03 & -&   1.7267e-03& - \\
$u_B=0.9$ & 240 & 8.7205e-04 &    2.4524   &  3.6879e-04 &   2.4529&  3.0632e-04&   2.4949\\
$\tau=1$ & 480&  1.2006e-04   &  2.8606   &  5.0480e-05   &  2.8690  &  4.1985e-05 &  2.8671\\
 & 960&  1.5942e-05  &   2.9129& 6.6663e-06 &    2.9208   &  5.1464e-06&   3.0282\\\hline
   & 120&  3.7573e-03   & -&  1.2122e-03 & -& 7.9211e-04 & - \\
$u_B=0.9$ & 240 & 7.4624e-04 &  2.3320 & 2.4164e-04  &    2.3267& 1.5061e-04 &  2.3949\\
$\tau=5$ & 480&  1.1994e-04   & 2.6373   &  3.8434e-05  &  2.6524  & 2.5089e-05  & 2.5857\\
 & 960&  1.5565e-05    & 2.9460 & 4.9190e-06  & 2.9660 &  3.1363e-06 &   2.9999\\\hline
   & 120&   2.1836e-03   &-& 9.1039e-04 &-&    5.7219e-04  &-\\
$u_B=\alpha$ & 240 &    3.2729e-04   & 2.7381&    1.3760e-04 &   2.7260  &8.9550e-05  &2.6757\\
$\tau=0.2$ & 480& 4.6856e-05  &  2.8043  &  1.9909e-05   & 2.7890  & 1.2935e-05  &2.7914\\
 & 960&   6.7382e-06  &  2.7978& 2.3182e-06    &3.1023 &  1.4109e-06 &  3.1965\\\hline
   & 120&   3.9014e-03  &-&  1.6388e-03   &-& 1.3873e-03 &-\\
$u_B=\alpha$ & 240 &  7.0517e-04  & 2.4680   & 2.9669e-04&   2.4656 &   2.4272e-04& 2.5149 \\
 $\tau=1$ & 480&  9.6528e-05   & 2.8690 & 4.0354e-05&   2.8781 &  3.3125e-05 & 2.8733\\
 & 960&    1.2890e-05 &   2.9047 &  5.3648e-06 &  2.9111   &  4.0754e-06&   3.0229\\\hline
    & 120&   3.0797e-03  &-&  9.9202e-04  &-&  6.4456e-04  &-\\
 $u_B=\alpha$ & 240 & 6.1133e-04  &   2.3328  &  1.9783e-04 &  2.3261&  1.2277e-04 &2.3924\\
 $\tau=5$ & 480& 9.7351e-05  & 2.6507 & 3.1222e-05   &2.6637   &  2.0263e-05 & 2.5990\\
& 960& 1.2396e-05  & 2.9733&  3.9513e-06 &   2.9822 &  2.4962e-06  &3.0210\\\hline
    & 120&  1.8244e-03  &-& 7.5548e-04 &-& 4.6671e-04  &-\\
 $u_B=0.75$ & 240 & 2.7262e-04&   2.7425  &  1.1419e-04   &  2.7260   & 7.3299e-05 &   2.6707 \\
$\tau=0.2$ & 480& 3.9198e-05 &  2.7980& 1.6562e-05  &   2.7855&  1.0681e-05 &   2.7788\\
 & 960& 5.4739e-06 &   2.8401  &1.9677e-06 &    3.0733  & 1.3232e-06&   3.0129\\\hline   
& 120&  3.2727e-03  &-& 1.3672e-03 &-&  1.1477e-03  &-\\
$u_B=0.75$ & 240 & 5.8671e-04  &  2.4798 & 2.4585e-04 &2.4754 & 1.9866e-04 &  2.5304\\
 $\tau=1$ & 480& 7.9974e-05 &  2.8750   & 3.3285e-05   &2.8848  &2.7033e-05 & 2.8775\\
 & 960&   1.0724e-05   & 2.8987 &   4.4466e-06 &  2.9041&  3.3341e-06 & 3.0193\\\hline
   & 120&   2.5902e-03 &-&  8.3335e-04  &-&  5.3882e-04 &-\\
$u_B=0.75$ & 240 & 5.1342e-04   & 2.3348   & 1.6611e-04&  2.3268& 1.0271e-04 &   2.3913\\
 $\tau=5$ & 480& 8.1062e-05 &   2.6630  &  2.6032e-05   & 2.6738   & 1.6813e-05 &   2.6109\\
& 960&    1.0173e-05  &  2.9944 & 3.2662e-06 &   2.9946&  2.0473e-06 & 3.0377\\\hline
  \end{tabular}
\end{footnotesize}
\end{center}
\caption{The accuracy test for the third order semi-discrete scheme for the MBL equation (\ref{MBL_new}) with 
$\epsilon=1$ and $M=2$.}
\label{3rd_order}
\end{table}

We will now use examples to study the solutions to MBL equation (\ref{MBL_new}) using the numerical schemes 
proposed 
in chapter \ref{central}.
We first notice that if we scale $t$ and $x$ as follows
\begin{eqnarray*}
\tilde{t}=\frac{t}{\epsilon},\qquad\tilde{x}=\frac{x}{\epsilon},
\end{eqnarray*}
then MBL (\ref{MBL_new}) equation can be written in terms of $\tilde{t}$ and $\tilde{x}$ as follows
\begin{eqnarray}
u_{\tilde{t}}+(f(u))_{\tilde{x}}=u_{\tilde{x}\tilde{x}}+\tau u_{\tilde{x}\tilde{x}\tilde{t}}.\label{scaling}
\end{eqnarray}
The scaled equation (\ref{scaling}) shows that it is the magnitude of $\frac{t}{\epsilon}$ and $\frac{x}{\epsilon}$ that determine the asymptotic behavior, not $t$, $x$, neither $\epsilon$ alone (\cite{duijn}). 
In addition,  (\ref{scaling}) also shows that the dispersive parameter $\tau$ denotes the relative importance of the dispersive term $u_{xxt}$. The bigger $\tau$ is, the more dispersive effect (\ref{MBL_new}) equation has. This can be seen from the computational results to be shown later in this section.

\begin{figure}[htbp]
\begin{center}
\includegraphics[scale=0.6]{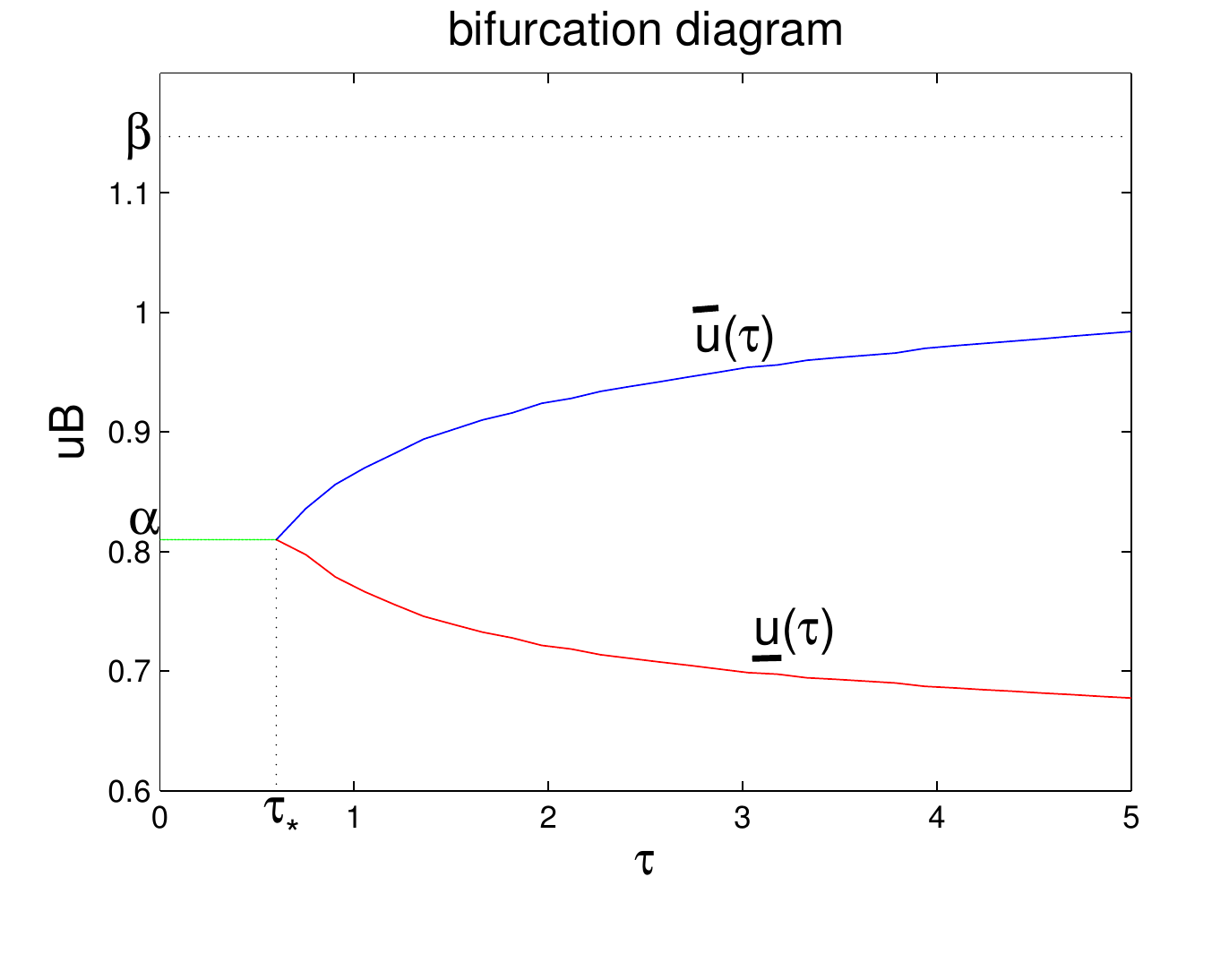}
\caption{The bifurcation diagram of the MBL equation (\ref{MBL}) with   the bifurcation parameters ($\tau, u_B$).}
\label{bif}
\end{center}
\end{figure}

\begin{figure}[htbp]
\begin{center}
\subfigure[]{\label{solution_type1}\includegraphics[width=0.325\textwidth]{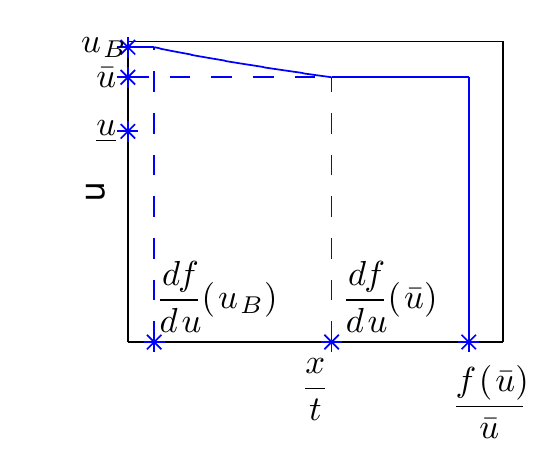}}
\subfigure[]{\label{solution_type2}\includegraphics[width=0.325\textwidth]{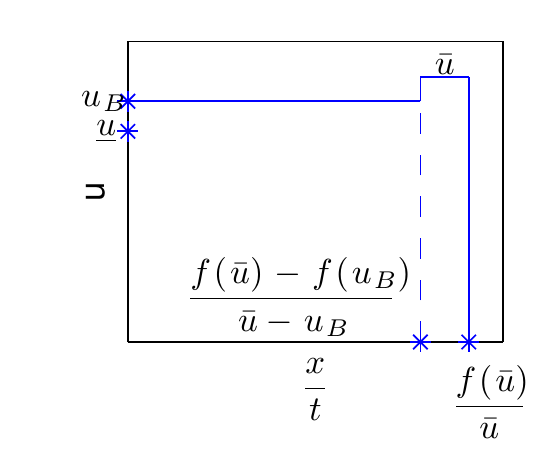}}
\subfigure[]{\label{solution_type3}\includegraphics[width=0.325\textwidth]{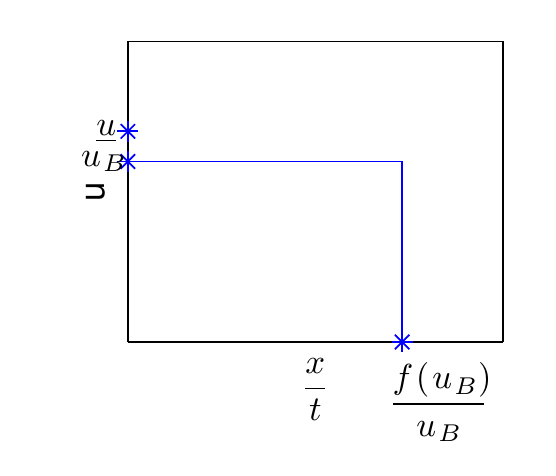}}
\caption{Given a fixed $\tau$, the three qualitatively different solution profiles due to different values of $u_B$. In particular, when $\tau>\tau_*$ and $\underline{u}<u_B<\bar{u}$, the solution profiles (Figure \ref{solution_type2}) displays non-monotonicity, which is consistent with the experimental observations (\cite{Dicarlo}). Figures \ref{solution_type1}, \ref{solution_type2} and \ref{solution_type3} are demonstrative figures. }
\label{solution_types}
\end{center}
\end{figure}

Duijn et al. \cite{duijn}  numerically provided a bifurcation diagram  (Figure \ref{bif}) of MBL (\ref{MBL_new}) equation as the dispersive parameter $\tau$ and the post-shock value $u_B$ of the initial condition vary.
The solution of (\ref{MBL_new}) has been proven to display qualitatively different profiles for parameter
values ($\tau, u_B$) falling in different regimes of the bifurcation diagram. 
In particular, for every fixed $\tau$ value, there are two critical $u_B$ values, namely, $\bar{u}$ and $\underline{u}$. From the bifurcation diagram (Figure \ref{bif}), it is clear that, when $\tau<\tau_*$, $\bar{u}=\underline{u}=\alpha$. For a fixed $\tau$ value, the solution has three different profiles.
\begin{itemize}
\item[(a)] If $u_B\in[\bar{u},1]$, the solution contains a plateau value $u_B$ for $0\le\frac{x}{t}\le\frac{df}{du}(u_B)$, a rarefaction wave connection $u_B$ to $\bar{u}$ for $\frac{df}{du}(u_B)\le\frac{x}{t}\le\frac{df}{du}(\bar{u})$, another plateau value $\bar{u}$ for $ \frac{df}{du}(\bar{u})<\frac{x}{t}<\frac{f(\bar{u})}{\bar{u}}$, and a shock from $\bar{u}$ down to $0$ at $\frac{x}{t}=\frac{f(\bar{u})}{\bar{u}}$ (see Figure \ref{solution_type1}).
\item[(b)]  If $u_B\in(\underline{u},\bar{u})$,  the solution contains a plateau value $u_B$ for $0\le\frac{x}{t}<\frac{f(\bar{u})-f(u_B)}{\bar{u}-u_B}$, a shock from $u_B$ up to $\bar{u}$ at $\frac{x}{t}=\frac{f(\bar{u})-f(u_B)}{\bar{u}-u_B}$, another plateau value $\bar{u}$ for $\frac{f(\bar{u})-f(u_B)}{\bar{u}-u_B}<\frac{x}{t}<\frac{f(\bar{u})}{\bar{u}}$, and a shock from $\bar{u}$ down to $0$ at $\frac{x}{t}=\frac{f(\bar{u})}{\bar{u}}$ (see Figure \ref{solution_type2}). The solution may exhibit a damped oscillation near $u=u_B$.
\item[(c)]  If $u_B\in(0,\underline{u}]$, the solution consists a single shock connecting $u_B$ and $0$ at $\frac{x}{t}=\frac{f(u_B)}{u_B}$ (see Figure \ref{solution_type3}). It may exhibit oscillatory behavior near $u=u_B$.
\end{itemize}
Notice that when $\tau>\tau_*$ and $\underline{u}<u_B<\bar{u}$, the solution profiles (\ref{solution_type2}) displays non-monotonicity, which is consistent with the experimental observations (\cite{Dicarlo}).

In the numerical computation we show below, we will therefore test the accuracy and capability of central schemes for different  parameter values ($\tau$ and $u_B$) that fall into various regimes of the bifurcation diagram, and
therefore display qualitatively different solution profiles.
The numerical experiments were carried out for $M=2$,  $\epsilon =0.001$ and $T=4000\times\epsilon$, i.e. $\tilde{T}=4000$ to get the asymptotic solution profiles, and $\Delta x$ was chosen to be $\frac{\epsilon}{10}$ and $\lambda=\frac{\Delta t}{\Delta x}$ was chosen to be 0.1.
The scheme used in the computation is the second order Trapezoid scheme as shown in section \ref{trapzoidsec}. The Midpoint scheme delivers similar computational results, hence is omitted here. The solution profiles at $\frac{T}{4}$ (blue), $\frac{2*T}{4}$ (green), $\frac{3*T}{4}$ (magenta) and $T$ (black) are chosen to demonstrate the time evolution of the solutions. The red dashed lines are used to denote the theoretical shock locations and plateau values for comparison purpose.

We start with $\tau>0$. Based on the bifurcation diagram (Figure \ref{bif}), we choose three representative $u_B$ values, i.e. $u_B=0.9>\alpha$, $u_B=\alpha=\sqrt{\frac{M}{M+1}}=\sqrt{\frac{2}{3}}$ (for $M=2$) and $u_B=0.75<\alpha$. For each fixed $u_B$, we choose three representative $\tau$ values, i.e. $\tau=0.2<\tau_*\approx0.61$, $\tau=1>\tau_*$  with $u_B=0.75<\underline{u}_{\tau=1}<u_B=\alpha<\bar{u}<u_B=0.9$, and $\tau=5$ with $u_B=0.75, \alpha, 0.9\in[\underline{u}_{\tau=5},\bar{u}_{\tau=5}]$.
We first use  this 9 pairs of $(\tau, u_B)$ values given in Table 
\ref{9examples} to validate the solution profiles with the demonstrative solution profiles given in Figure \ref{solution_types}.
\begin{table}[htbp]
\begin{center}
\begin{tabular}{|l|l|l|l|}
\hline
  $(\tau,u_B)$     & Example 4 & Example 5 & Example 6 \\ \hline
 Example 1         & $(0.2, 0.9)$& $(1, 0.9)$&$ ( 5, 0.9)$\\\hline
 Example 2         & $(0.2, \alpha)$& $ (1, \alpha)$& $(5, \alpha)$\\\hline
 Example 3         & $(0.2, 0.75)$&$(1, 0.75)$&$(5, 0.75)$\\\hline
\end{tabular}
\caption{9 pairs of $(\tau, u_B)$ values with either fixed $\tau$ value or 
fixed $u_B$ value used in Examples 1 -- 6.}
\label{9examples}
\end{center}
\end{table}

\begin{figure}[htbp]
\subfiguretopcaptrue
\subfigure[$(\tau, 
u_B)=(0.2, 
0.9)$]{\label{trapezoid_11_epsilon_0.001}\includegraphics[width=0.325\textwidth]{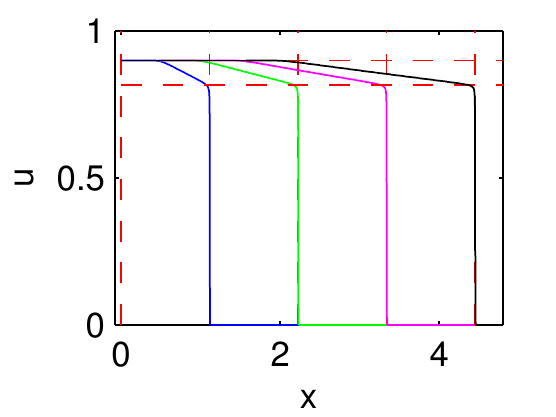}}
\subfigure[$(\tau, 
u_B)=(1, 
0.9)$]{\label{trapezoid_12_epsilon_0.001}\includegraphics[width=0.325\textwidth]{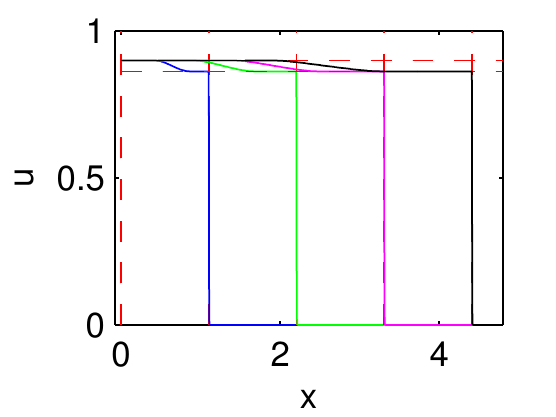}}
\subfigure[$(\tau, 
u_B)=(5, 
0.9)$]{\label{trapezoid_13_epsilon_0.001}\includegraphics[width=0.325\textwidth]{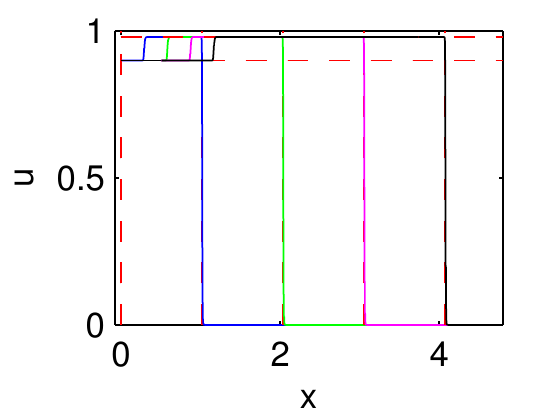}}\\
\subfigure[$(\tau, 
u_B)=(0.2, 
\alpha)$]{\label{trapezoid_21_epsilon_0.001}\includegraphics[width=0.325\textwidth]{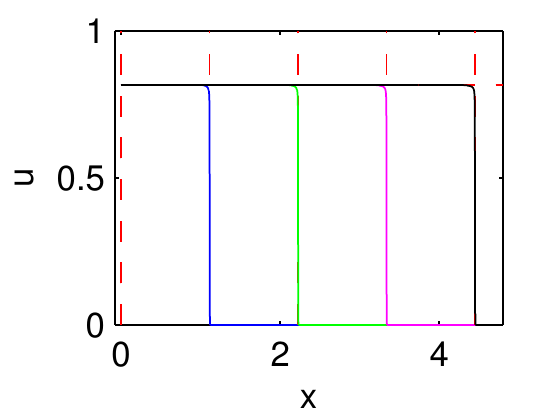}}
\subfigure[$(\tau, 
u_B)=(1, 
\alpha)$]{\label{trapezoid_22_epsilon_0.001}\includegraphics[width=0.325\textwidth]{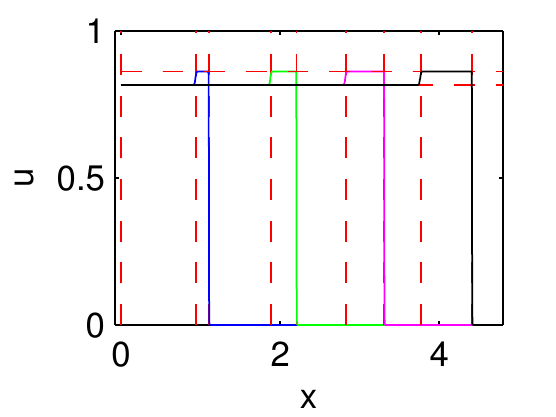}}
\subfigure[$(\tau, 
u_B)=(5, 
\alpha)$]{\label{trapezoid_23_epsilon_0.001}\includegraphics[width=0.325\textwidth]{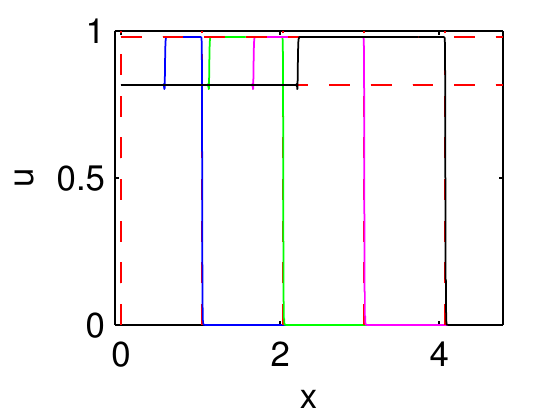}}\\
\subfigure[$(\tau, 
u_B)=(0.2, 
0.75)$]{\label{trapezoid_31_epsilon_0.001}\includegraphics[width=0.325\textwidth]{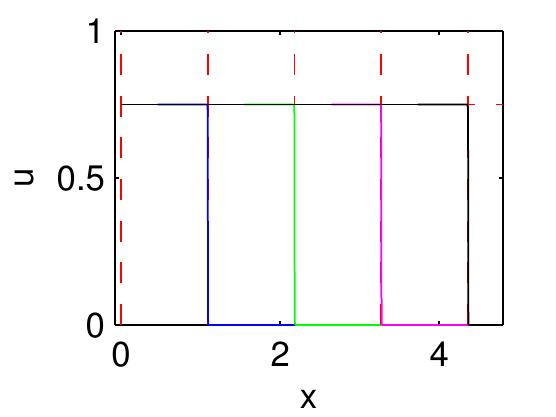}}
\subfigure[$(\tau, 
u_B)=(1, 
0.75)$]{\label{trapezoid_32_epsilon_0.001}\includegraphics[width=0.325\textwidth]{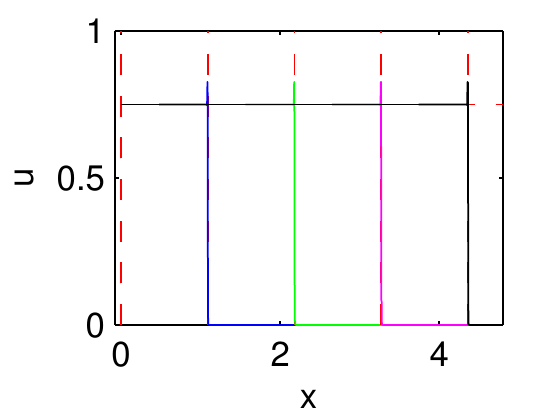}}
\subfigure[$(\tau, 
u_B)=(5, 
0.75)$]{\label{trapezoid_33_epsilon_0.001}\includegraphics[width=0.325\textwidth]{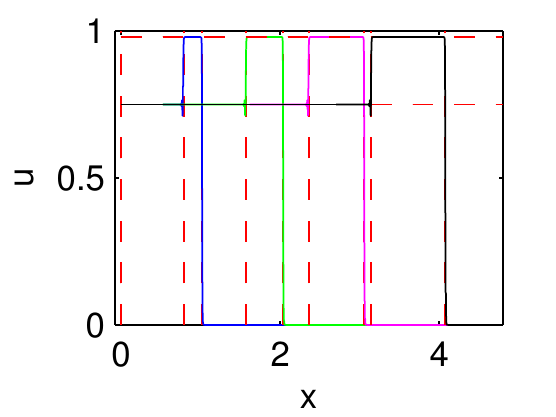}}
\caption{Numerical solutions to MBL equation with parameter settings  fall in different regimes of the 
bifurcation diagram (Figure \ref{bif}). 
The color coding is for different time: $\frac{1}{4}T$ (blue), $\frac{2}{4}T$ (green), $\frac{3}{4}T$ (magenta) 
and $T$ (black).
The results are discussed in examples 1 -- 6. In figures \ref{trapezoid_21_epsilon_0.001} -- 
\ref{trapezoid_23_epsilon_0.001}, $\alpha = \sqrt{\frac{M}{M+1}} = \sqrt{\frac{2}{3}}$ for $M=2$.}
 \label{nine} 
\end{figure}

\noindent{\bf Example 1} $(\tau, u_B) = (0.2, 0.9),(\tau, u_B) = (1, 0.9),(\tau, u_B) = (5, 0.9)$.\\
When $u_B=0.9>\alpha$ is fixed, we increase $\tau$ from 0.2 to 1 to 5 (Figure \ref{trapezoid_11_epsilon_0.001} , \ref{trapezoid_12_epsilon_0.001} , \ref{trapezoid_13_epsilon_0.001}), the dispersive effect starts to dominate the solution profile. When $\tau = 0.2$ (Figure \ref{trapezoid_11_epsilon_0.001}), the solution profile is similar to the classical BL equation solution (see Figure \ref{BL2}), with a rarefaction wave for $\frac{x}{t}\in[f'(u=0.9),f'(u=\alpha)=f'(u=\bar{u}_{\tau=0.2})]$ and a shock from $u=\alpha$ to $u=0$ at $\frac{x}{t}=f'(\alpha)$. 
This corresponds to Figure \ref{solution_type1} with $\frac{df}{du}(\bar{u}_{\tau=0.2}=\alpha)=\frac{f(\bar{u}_{\tau=0.2})}{\bar{u}_{\tau =0.2}}=\frac{f(\alpha)}{\alpha}$.
When $\tau = 1$ (Figure \ref{trapezoid_12_epsilon_0.001}), the rarefaction wave is between $\frac{x}{t}\in[f'(u=0.9), f'(u=\bar{u}_{\tau=1})]$ and the solution remains at the plateau value $u=\bar{u}_{\tau=1}$ for $\frac{x}{t}\in[ f'(u=\bar{u}_{\tau=1}),\frac{f(\bar{u}_{\tau=1})}{\bar{u}_{\tau=1}}]$ and the shock occurs at $\frac{x}{t}=\frac{f(\bar{u}_{\tau=1})}{\bar{u}_{\tau=1}}$. 
This corresponds to Figure \ref{solution_type1} with $u_B=0.9>\bar{u}_{\tau=1}\approx 0.86$.
When $\tau =5$ (Figure \ref{trapezoid_13_epsilon_0.001}), the solution displays the first shock from $u=0.9$ to $u=\bar{u}_{\tau=5}$ at $\frac{x}{t}=\frac{f(\bar{u}_{\tau=5})-f(u_B)}{\bar{u}_{\tau=5}-u_B}$, and then remains at the plateau value $u=\bar{u}_{\tau=5}$ for $\frac{x}{t}\in[ \frac{f(\bar{u}_{\tau=5})-f(u_B)}{\bar{u}_{\tau=5}-u_B},\frac{f(\bar{u}_{\tau=5)}}{\bar{u}_{\tau=5}}]$ and the second shocks occurs at $\frac{x}{t}=\frac{f(\bar{u}_{\tau=5)}}{\bar{u}_{\tau=5}}$.
This corresponds to Figure \ref{solution_type2} with $\underline{u}_{\tau=5}\approx 0.68<u_B=0.9<\bar{u}_{\tau=5}\approx 0.98$. Notice that as $\tau$ increases, the rarefaction region shrinks and the plateau region enlarges.\\

\noindent{\bf Example 2} $(\tau, u_B) = (0.2, \alpha),(\tau, u_B) = (1,\alpha),(\tau, u_B) = (5, \alpha)$.\\
When $u_B=\alpha$ is fixed, we increase $\tau$ from 0.2 to 1 to 5 (Figure \ref{trapezoid_21_epsilon_0.001} , \ref{trapezoid_22_epsilon_0.001} , \ref{trapezoid_23_epsilon_0.001}), the dispersive effect starts to dominate the solution profile.  When $\tau=0.2$, the solution displays one single shock at $\frac{x}{t}=\frac{f(\alpha)}{\alpha}$. For both $\tau=1$ and $\tau=5$, the solution has two shocks, one at $\frac{x}{t}=\frac{f(\bar{u}_{\tau=1(\tau=5 \mathrm{~ respectively})})-f(\alpha)}{\bar{u}_{\tau=1(\tau=5 \mathrm{~respectively})}-\alpha}$, and another one at $\frac{x}{t}=\frac{f(\bar{u}_{\tau=1(\tau=5\mathrm{~ respectively})})}{\bar{u}_{\tau=1(\tau=5 \mathrm{~respectively})}}$. For both $\tau =1$ and $\tau =5$ (Figures \ref{trapezoid_22_epsilon_0.001} \ref{trapezoid_23_epsilon_0.001}), the solutions correspond to Figure \ref{solution_type2}, which are consistent to the experimental observations. Notice that as $\tau$ increases from 1 to 5, i.e., the dispersive effect increases, the inter-shock interval length increases at every fixed time (compare Figure \ref{trapezoid_22_epsilon_0.001}  with Figure \ref{trapezoid_23_epsilon_0.001}). In addition, for  fix $\tau=1$ ($\tau=5$ respectively), as time progresses, the inter-shock interval length increases in the linear fashion (see Figure \ref{trapezoid_22_epsilon_0.001} (Figure \ref{trapezoid_23_epsilon_0.001} respectively) ).\\

\noindent{\bf Example 3} $(\tau, u_B) = (0.2, 0.75),(\tau, u_B) = (1,0.75),(\tau, u_B) = (5, 0.75)$.\\
When $u_B=0.75<=\alpha$ is fixed, we increase $\tau$ from 0.2 to 1 to 5 (Figure \ref{trapezoid_31_epsilon_0.001} , \ref{trapezoid_32_epsilon_0.001} , \ref{trapezoid_33_epsilon_0.001}), the dispersive effects starts to dominate the solution profile in the similar fashion as $u_B=0.9$ and $u_B=\alpha$. Notice that when $\tau=1$, since $u_B=0.75$ is very close to $\underline{u}_{\tau=1}$, the solution displays oscillation at $\frac{x}{t}=\frac{f(u_B)}{u_B}$ (Figure  \ref{trapezoid_32_epsilon_0.001}). If we increase $\tau$ further to $\tau=5$, the dispersive effect is strong enough to create a plateau value at $\bar{u}\approx 0.98$ (see Figure \ref{trapezoid_33_epsilon_0.001}).\\

\noindent{\bf Example 4} $(\tau, u_B) = (0.2, 0.9),(\tau, u_B) = (0.2,\alpha),(\tau, u_B) = (0.2, 0.75)$.\\
Now, we fix $\tau=0.2$, decrease $u_B$ from 0.9 to $\alpha$, to 0.75 
(Figures\ref{trapezoid_11_epsilon_0.001} \ref{trapezoid_21_epsilon_0.001} \ref{trapezoid_31_epsilon_0.001}). If $u_B>\alpha$ the solution consists a rarefaction wave connecting $u_B$ down to $\alpha$, then a shock from $\alpha$ to 0, otherwise, the solution consists a single shock from $u_B$ down to 0. In all cases, since $\tau = 0.2<\tau_*$, regardless of the $u_B$ value, the solution will not display non-monotone behavior, due to the lack of dispersive effect.\\

\noindent{\bf Example 5} $(\tau, u_B) = (1, 0.9),(\tau, u_B) = (1,\alpha),(\tau, u_B) = (1, 0.75)$.\\
Now, we fix $\tau=1$, decrease $u_B$ from 0.9 to $\alpha$, to 0.75 (Figures\ref{trapezoid_12_epsilon_0.001} \ref{trapezoid_22_epsilon_0.001} \ref{trapezoid_32_epsilon_0.001}).  If $u_B=0.9>\bar{u}_{\tau=1}$, the solution consists a rarefaction wave connecting $u_B$ and $\bar{u}$, and a shock connecting $\bar{u}$ down to 0 (Figure \ref{trapezoid_12_epsilon_0.001}). Even if $\underline{u}<u_B<\bar{u}$, because $\tau=1>\tau_*$, the solution still has a chance to increase to the plateau value $\bar{u}$ as seen in Figure \ref{trapezoid_22_epsilon_0.001}. But, if $u_B$ is too small, for example, $u_B=0.75<\underline{u}$, the solution does not increase to $\bar{u}$ any more, instead, it consists a single shock connecting $u_B$ down to 0 (Figure \ref{trapezoid_32_epsilon_0.001}). \\

\noindent{\bf Example 6} $(\tau, u_B) = (5, 0.9),(\tau, u_B) = (5,\alpha),(\tau, u_B) = (5, 0.75)$.\\
Now, we fix $\tau=5$, decrease $u_B$ from 0.9 to $\alpha$, to 0.75 (Figures\ref{trapezoid_13_epsilon_0.001} \ref{trapezoid_23_epsilon_0.001} \ref{trapezoid_33_epsilon_0.001}).
For all three $u_B$, they are between $\underline{u}_{\tau=5}$ and $\bar{u}_{\tau=5}$, hence all increase to the plateau value $\bar{u}_{\tau=5}\approx 0.98$ before dropping to 0. Notice that as $u_B$ decreases, the inter-shock interval length decreases at every fixed time (compare Figures \ref{trapezoid_13_epsilon_0.001}, \ref{trapezoid_23_epsilon_0.001} and \ref{trapezoid_33_epsilon_0.001}). This shows that when the dispersive effect is strong ($\tau>\tau_*$), the bigger $u_B$ is, the bigger region the solution stays at the plateau value.
\\

\begin{figure}[htbp]
\subfiguretopcaptrue
\subfigure[$(\tau, u_B)=(0, 0.9)$]{%
\label{trapezoid_61_epsilon_0.001}%
\includegraphics[width=0.325\textwidth]{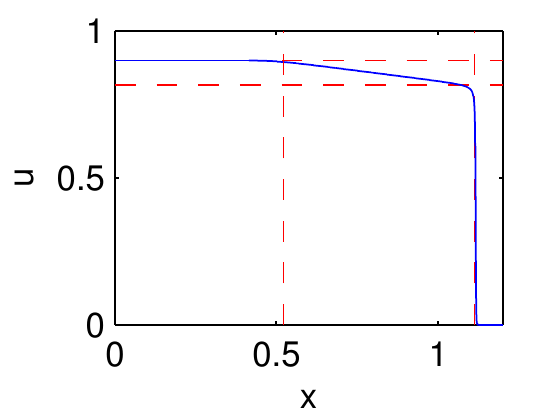}}
\subfigure[$(\tau, u_B)=(0, \alpha)$]{%
\label{trapezoid_62_epsilon_0.001}%
\includegraphics[width=0.325\textwidth]{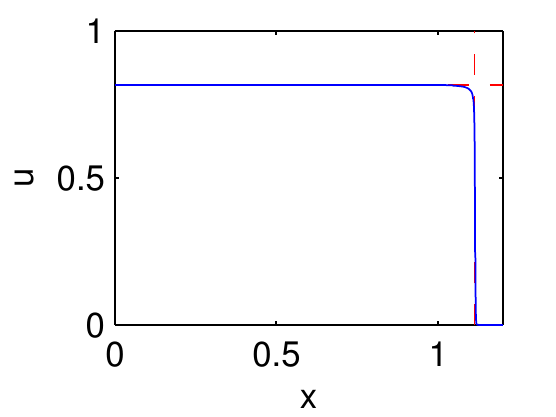}}
\subfigure[$(\tau, u_B)=(0, 0.75)$]{%
\label{trapezoid_63_epsilon_0.001}%
\includegraphics[width=0.325\textwidth]{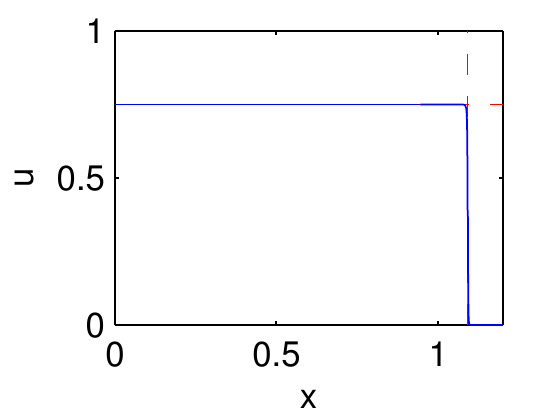}}
\caption{The numerical solutions of the MBL equation at T = 1 with $\tau=0$ and different $u_B$ values. The 
results are discussed in example 7. 
}
\label{trapezoid_6}
\end{figure}

\noindent{\bf Example 7} $(\tau, u_B) = (0, 0.9), (\tau, u_B) = (0, \alpha),  (\tau, u_B) = (0, 0.75)$.\\
We now show the solution profiles for the extreme $\tau$ value, i.e. $\tau =0$ in Figures  
\ref{trapezoid_61_epsilon_0.001} ($u_B = 0.9$),  \ref{trapezoid_62_epsilon_0.001} ($u_B = 
\alpha$) and \ref{trapezoid_63_epsilon_0.001} ($u_B = 0.75$). Notice that these are cases of 
classical BL equation with small diffusion $\epsilon u_{xx}$. We compare Figures 
\ref{trapezoid_61_epsilon_0.001}, \ref{trapezoid_62_epsilon_0.001} and \ref{trapezoid_63_epsilon_0.001} with the 
solution of the classical BL equation  given in Figures \ref{BL1} and \ref{BL2}, it is clear that they show 
qualitatively same solution profiles. The difference is that due to the 
diffusion term in the MBL equation, as shown in Figure \ref{trapezoid_6}, the 
solutions do not have sharp edges right at the shock, instead, the solutions 
smear out a little.
If we compare Figures \ref{trapezoid_61_epsilon_0.001}, \ref{trapezoid_62_epsilon_0.001} and 
\ref{trapezoid_63_epsilon_0.001} with Figures \ref{trapezoid_11_epsilon_0.001}, 
\ref{trapezoid_21_epsilon_0.001} and \ref{trapezoid_31_epsilon_0.001}, there is no visible difference. This shows 
that once $\tau<\tau_*$, solution profile will stay the same for a fixed $u_B$ value.\\

\begin{figure}[htbp]
\subfiguretopcaptrue
\subfigure[$(\tau, u_B)=(5, 0.99)$]{%
\label{trapezoid_41_epsilon_0.001}%
\includegraphics[width=0.325\textwidth]{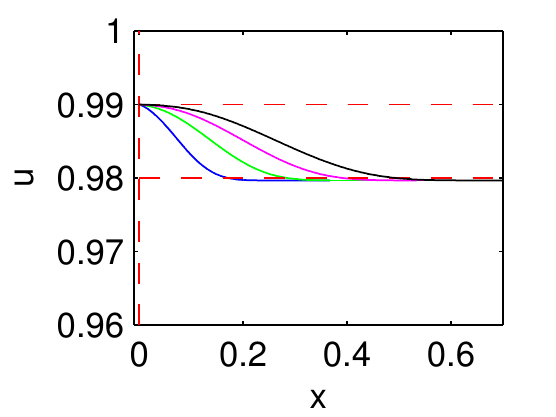}}
\subfigure[$(\tau, u_B)=(5, 0.98)$]{%
\label{trapezoid_42_epsilon_0.001}%
\includegraphics[width=0.325\textwidth]{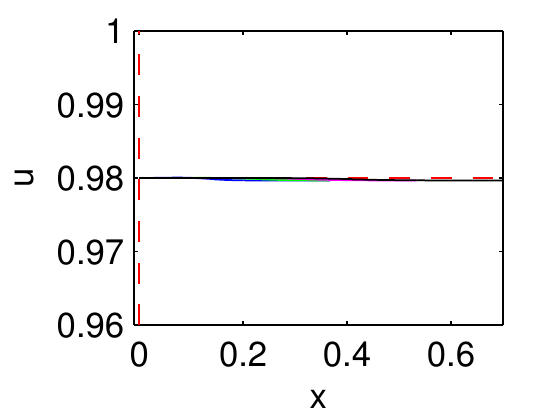}}
\subfigure[$(\tau, u_B)=(5, 0.97)$]{%
\label{trapezoid_43_epsilon_0.001}%
\includegraphics[width=0.325\textwidth]{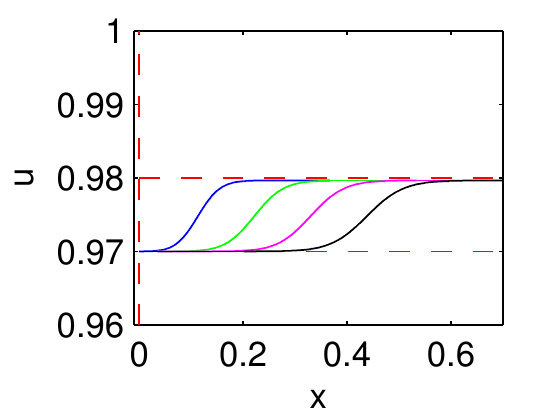}}
\caption{Numerical solutions to MBL equation with $u_B$ close to $\bar{u}_{\tau=5}\approx 0.98$.
The color coding is for different time: $\frac{1}{4}T$ (blue), $\frac{2}{4}T$ (green), $\frac{3}{4}T$ (magenta) 
and $T$ (black). The results are discussed in example 8.}

\label{nine_41_42_43}
\end{figure}

\noindent{\bf Example 8} $(\tau, u_B) = (5, 0.99),(\tau, u_B) = (5,0.98),(\tau, u_B) = (5, 0.97)$.\\
We also study the solution profiles for $u_B$ close to $\bar{u}$. For example, when $\tau=5$, $\bar{u}\approx 0.98$, we hence choose $u_B=0.99$, $u_B=0.98$, $u_B=0.97$ and solutions are shown in Figure \ref{trapezoid_41_epsilon_0.001}, \ref{trapezoid_42_epsilon_0.001}, \ref{trapezoid_43_epsilon_0.001}.
If $u_B=0.99>\bar{u}_{\tau=5}\approx 0.98$, the solution drops to the plateau value $\bar{u}$, then drops to 0 (see Figure \ref{trapezoid_41_epsilon_0.001}).  If $u_B=0.98\approx \bar{u}_{\tau=5}$, the solution remains at plateau value $\bar{u}_{\tau=5}$ and then drop to 0 (see Figure \ref{trapezoid_42_epsilon_0.001}). If $u_B=0.97<\bar{u}_{\tau=5}$, the solution increases to the plateau value $\bar{u}_{\tau=5}\approx 0.98$, then drops to 0. In all cases, the transition from $u_B$ to $\bar{u}_{\tau=5}\approx 0.98$ takes very small space. In the majority space, the solution keeps to be the plateau value $\bar{u}_{\tau=5}\approx 0.98$. \\

\noindent{\bf Example 9} $(\tau, u_B) = (5, 0.7)$, 
$(\tau, u_B) = (5,0.69)$,
$(\tau, u_B) = (5, 0.68)$,
$(\tau, u_B) = (5,0.67)$,
$(\tau, u_B) = (5,0.66)$.

\begin{figure}[hptb]
\subfiguretopcaptrue
\subfigure[$(\tau, u_B)=(5, 0.7)$]{%
\label{trapezoid_51_epsilon_0.001}%
\includegraphics[width=0.325\textwidth]{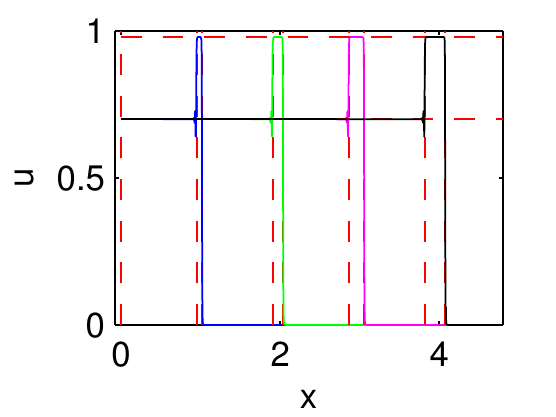}}
\subfigure[$(\tau, u_B)=(5, 0.69)$]{%
\label{trapezoid_52_epsilon_0.001}%
\includegraphics[width=0.325\textwidth]{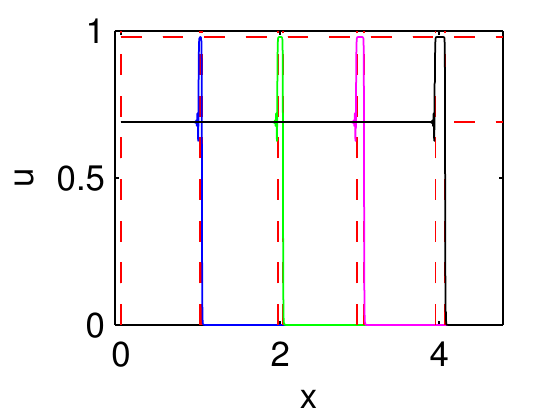}}
\subfigure[$(\tau, u_B)=(5, 0.68)$]{%
\label{trapezoid_53_epsilon_0.001}%
\includegraphics[width=0.325\textwidth]{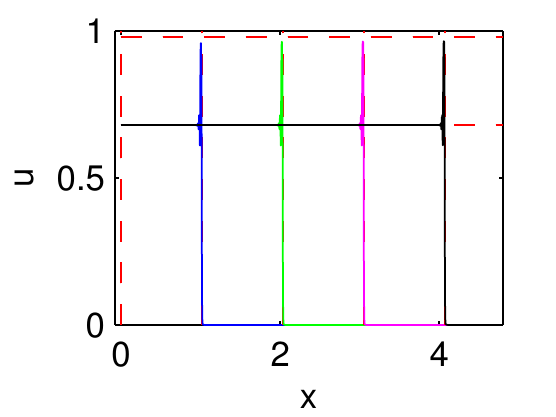}}
\subfigure[$(\tau, u_B)=(5, 0.67)$]{%
\label{trapezoid_54_epsilon_0.001}%
\includegraphics[width=0.325\textwidth]{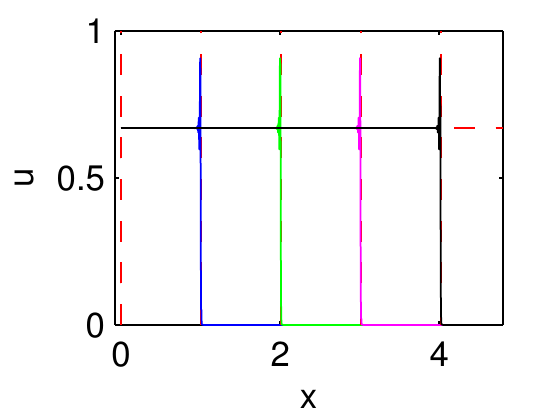}}
\subfigure[$(\tau, u_B)=(5, 0.66)$]{%
\label{trapezoid_55_epsilon_0.001}%
\includegraphics[width=0.325\textwidth]{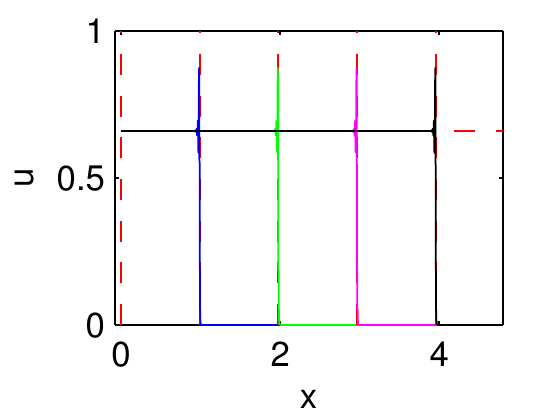}}
\caption{Numerical solutions to MBL equation with $u_B$ close to 
\underline{$u$}$_{\tau=5}\approx 0.68$. 
The color coding is for different time: $\frac{1}{4}T$ (blue), $\frac{2}{4}T$ (green), 
$\frac{3}{4}T$ (magenta) and $T$ (black). The results are discussed in example 9.}
\label{nine_51_52_53_54_55} 
\end{figure}

\noindent In addition, we  study the solution profiles for $u_B$ close to $\underline{u}$. For example, when 
$\tau=5$, $\underline{u}\approx 0.68$, we hence choose $u_B=0.7$, $u_B=0.69$, $u_B=0.68$, $u_B=0.67$, $u_B=0.66$ and solutions are shown in Figures \ref{trapezoid_51_epsilon_0.001},  \ref{trapezoid_52_epsilon_0.001}, \ref{trapezoid_53_epsilon_0.001},  \ref{trapezoid_54_epsilon_0.001}, \ref{trapezoid_55_epsilon_0.001}. As $u_B$ decreases crossing $\underline{u}_{\tau=5}\approx 0.68$, the solution gradually stops increasing to the plateau value $\bar{u}_{\tau=5}$, and the inter-shock interval length decreases (compare Figures \ref{trapezoid_51_epsilon_0.001}, \ref{trapezoid_52_epsilon_0.001} and \ref{trapezoid_53_epsilon_0.001}). The oscillation in Figures \ref{trapezoid_54_epsilon_0.001} and \ref{trapezoid_55_epsilon_0.001} are due to the fact that $u_B$ values are too close to $\underline{u}_{\tau=5}$. This confirms that even with big dispersive effect (say $\tau=5$), if $u_B$ is too small (e.g. $u_B<\underline{u}$), the solution will not exhibit non-monotone behavior.\\

\begin{figure}[!h]
\subfiguretopcaptrue
\subfigure[$(\tau, u_B)=(0.2, 
0.6)$]{\label{oscillation_71}\includegraphics[width=0.325\textwidth]{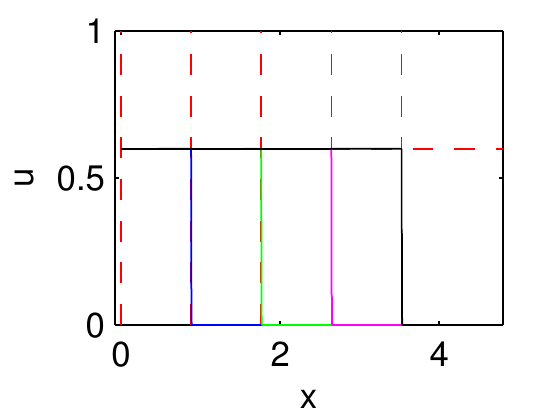}}
\subfigure[$(\tau, u_B)=(1, 
0.6)$]{\label{oscillation_72}\includegraphics[width=0.325\textwidth]{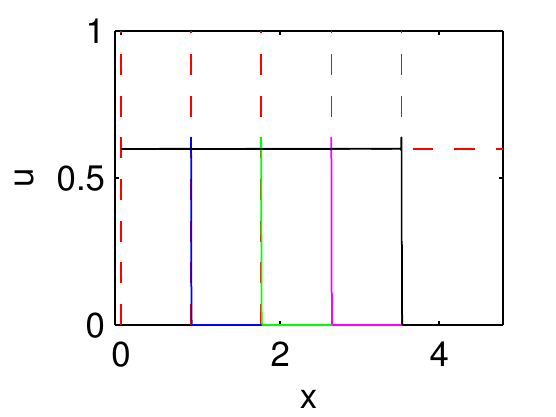}}
\subfigure[ $(\tau, u_B)=(5, 
0.6)$]{\label{oscillation_73}\includegraphics[width=0.325\textwidth]{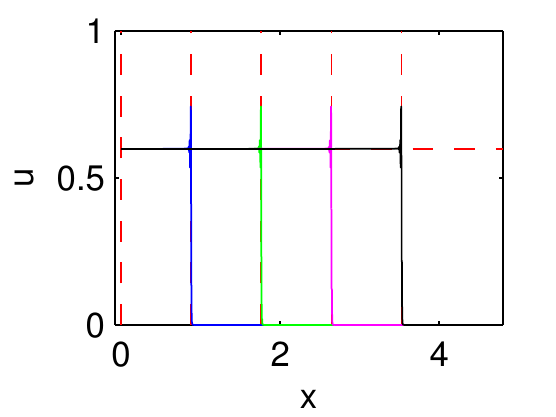}}
\subfigure[Fig \ref{oscillation_71} zoomed in at 
$\frac{1}{4}T$]{\label{oscillation_71_zoom1}\includegraphics[width=0.325\textwidth]{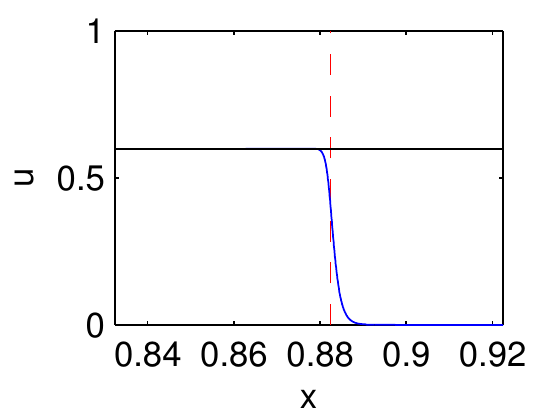}}
\subfigure[Fig \ref{oscillation_72} zoomed in at 
$\frac{1}{4}T$]{\label{oscillation_72_zoom1}\includegraphics[width=0.325\textwidth]{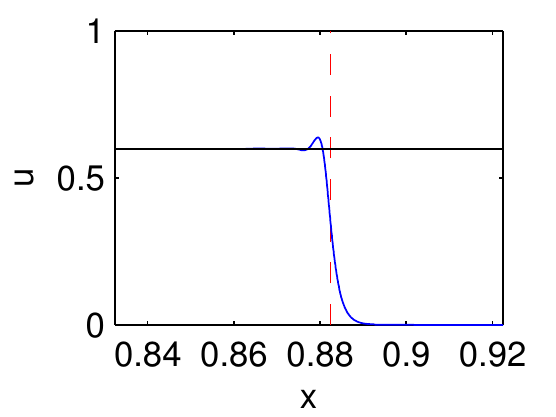}}
\subfigure[Fig \ref{oscillation_73} zoomed in at 
$\frac{1}{4}T$]{\label{oscillation_73_zoom1}\includegraphics[width=0.325\textwidth]{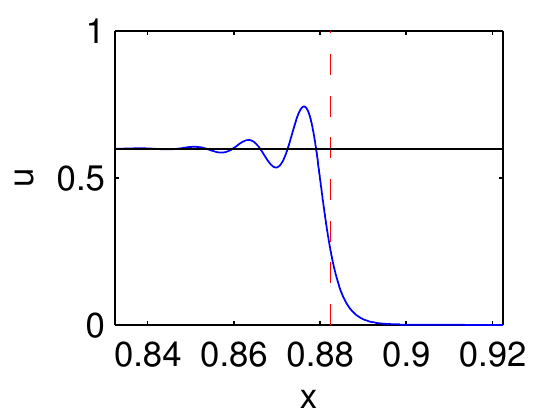}}
\subfigure[Fig \ref{oscillation_71} zoomed in at 
$T$]{\label{oscillation_71_zoom4}\includegraphics[width=0.325\textwidth]{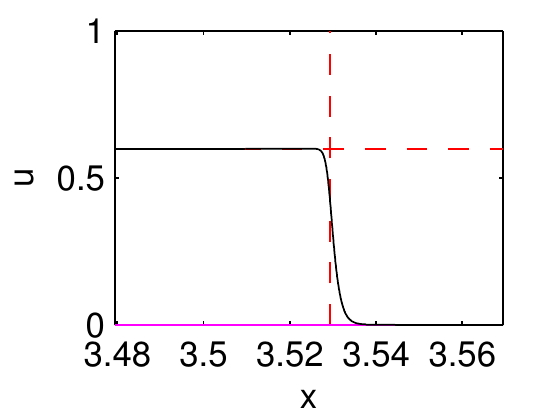}}
\subfigure[Fig \ref{oscillation_72} zoomed in at 
$T$]{\label{oscillation_72_zoom4}\includegraphics[width=0.325\textwidth]{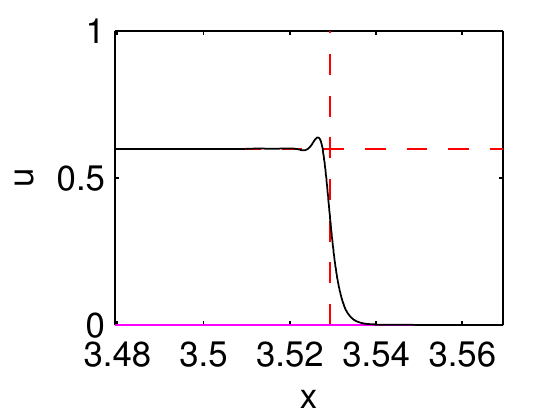}}
\subfigure[Fig \ref{oscillation_73} zoomed in at 
$T$]{\label{oscillation_73_zoom4}\includegraphics[width=0.325\textwidth]{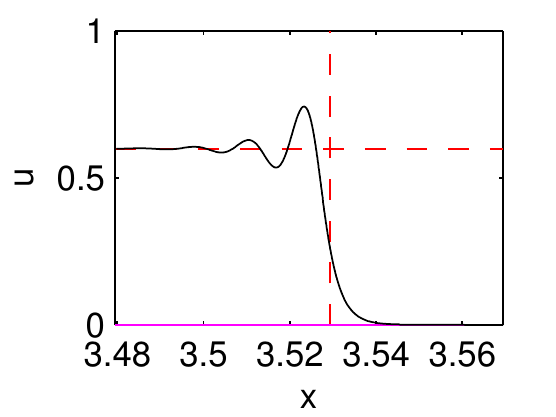}}
\caption{Numerical solutions to MBL equation with  small constant $u_B=0.6$ and different $\tau$ values.
The figures on the second and third rows are the magnified versions of the first row at $t=\frac{1}{4}T$ and 
$t=T$ respectively.
The color coding is for different time: $\frac{1}{4}T$ (blue), $\frac{2}{4}T$ (green), $\frac{3}{4}T$ (magenta) 
and $T$ (black).
The results are discussed in examples 10. 
}
\label{oscillation}
\end{figure}

\noindent{\bf Example 10} $(\tau, u_B) = (0.2, 0.6)$, $(\tau, u_B) = (1,0.6)$, $(\tau, u_B) = (5, 0.6)$.\\
We fix $u_B$ to be small, and in this example, we take it to be $u_B=0.6$. We vary the $\tau$ value, from $\tau = 
0.2<\tau_*$ to $\tau = 1$  barely larger than $\tau_*$ to $\tau = 5>\tau_*$. The numerical solutions are given in 
Figure \ref{oscillation_71}, \ref{oscillation_72}, \ref{oscillation_73}.  As $\tau$ increases, the post-shock 
value remains the same, but there will be oscillation generated as $\tau$ becomes larger than $\tau_*$. Figures 
\ref{oscillation_71_zoom1}, \ref{oscillation_72_zoom1} and \ref{oscillation_73_zoom1} show that as $\tau$ 
increases, the oscillation amplitude increases and oscillates more rounds. 
Notice that $\tau$ is the dispersive parameter, and this means that even for small $u_B$ value, different 
dispersive parameter values still give different dispersive effects, although none can bring the solution to the 
plateau value $\bar{u}$.
Comparing 
Figures \ref{oscillation_71_zoom1}, \ref{oscillation_72_zoom1} and \ref{oscillation_73_zoom1}  with
Figures \ref{oscillation_71_zoom4}, \ref{oscillation_72_zoom4} and \ref{oscillation_73_zoom4}, it is clear that 
the oscillation amplitude remains steady with respect to time.\\

\begin{figure}[hptb]
\subfiguretopcaptrue
\subfigure[$(\tau, 
u_B)=(0.2, 
0.9)$]{\label{diff_epsilon_11}\includegraphics[width=0.325\textwidth]{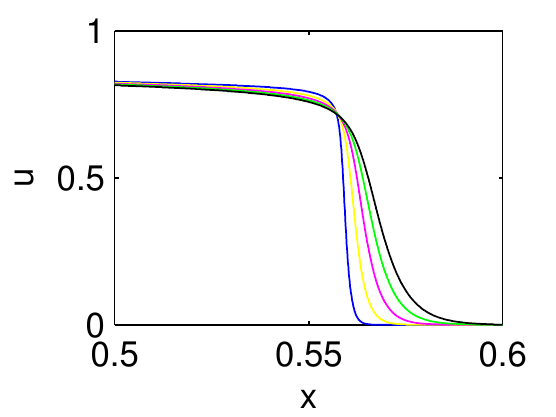}}
\subfigure[$(\tau, 
u_B)=(1, 
0.9)$]{\label{diff_epsilon_12}\includegraphics[width=0.325\textwidth]{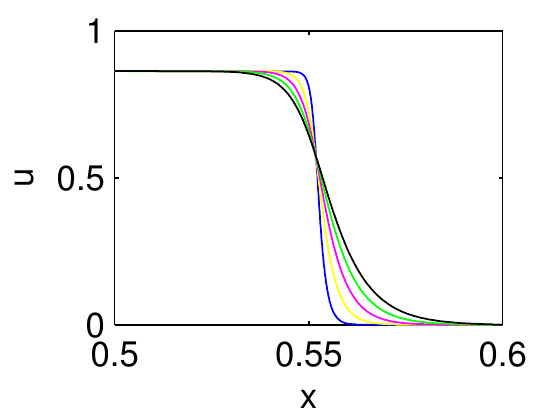}}
\subfigure[$(\tau, 
u_B)=(5, 
0.9)$]{\label{diff_epsilon_13}\includegraphics[width=0.325\textwidth]{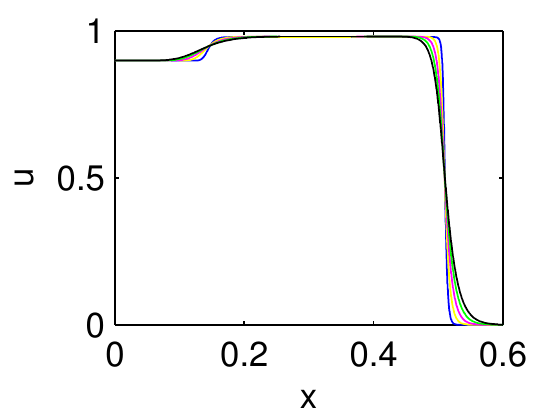}}
\subfigure[ 
$(\tau, 
u_B)=(0.2, 
\alpha)$]{\label{diff_epsilon_21}\includegraphics[width=0.325\textwidth]{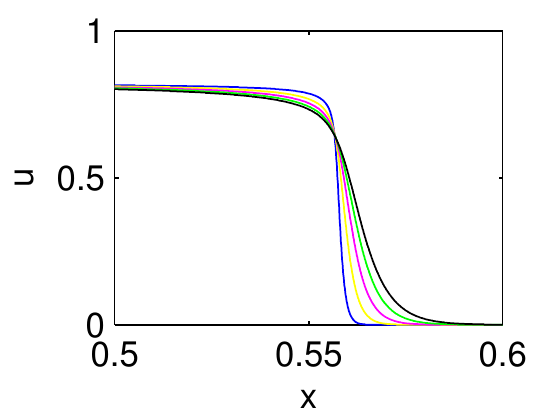}}
\subfigure[$(\tau, 
u_B)=(1, 
\alpha)$]{\label{diff_epsilon_22}\includegraphics[width=0.325\textwidth]{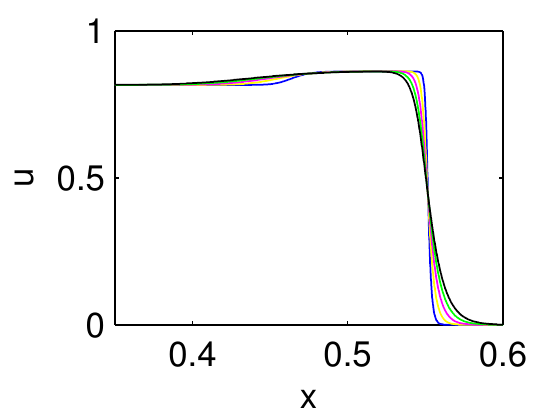}}
\subfigure[$(\tau, 
u_B)=(5, 
\alpha)$]{\label{diff_epsilon_23}\includegraphics[width=0.325\textwidth]{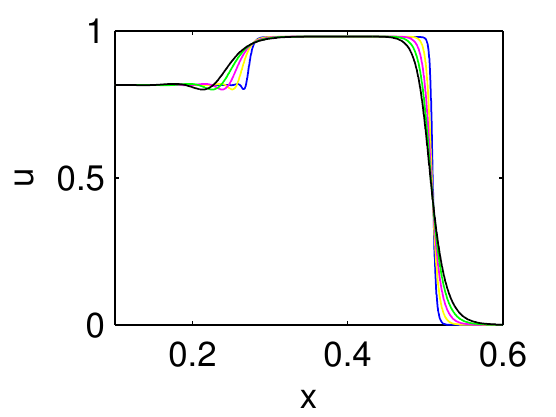}}
\subfigure[$(\tau, 
u_B)=(0.2, 
0.75)$]{\label{diff_epsilon_31}\includegraphics[width=0.325\textwidth]{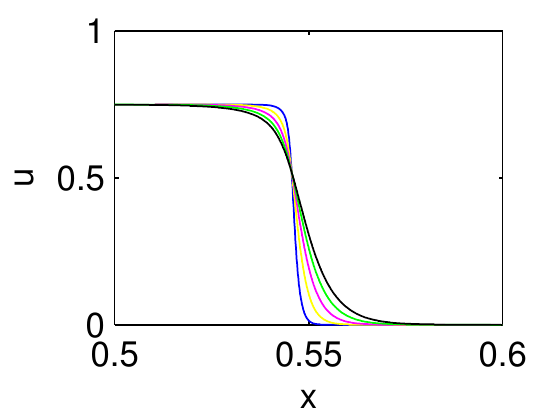}}
\subfigure[$(\tau, 
u_B)=(1, 
0.75)$]{\label{diff_epsilon_32}\includegraphics[width=0.325\textwidth]{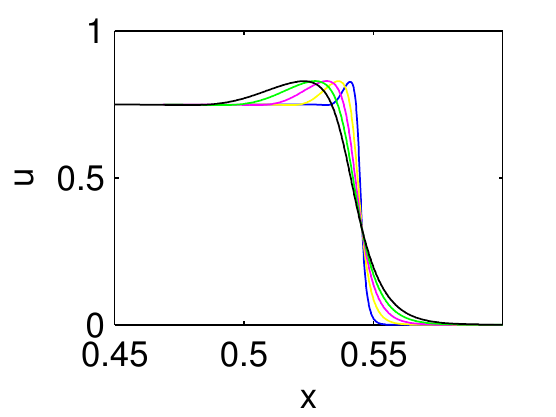}}
\subfigure[$(\tau, 
u_B)=(5, 
0.75)$]{\label{diff_epsilon_33}\includegraphics[width=0.325\textwidth]{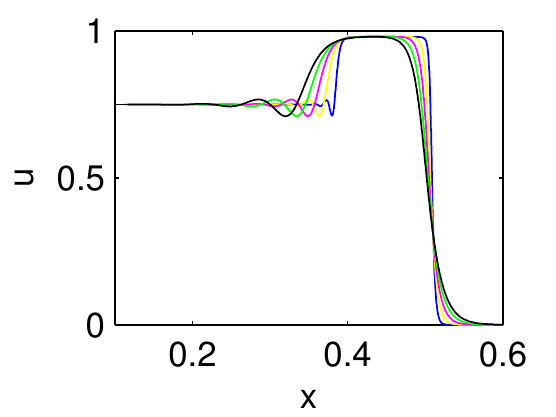}}
\caption{The numerical solutions of MBL equation at $T=0.5$ with $\epsilon= 0.001$ (blue), 
$\epsilon=0.002$ (yellow), $\epsilon=0.003$ (magenta),  $\epsilon=0.004$ (green), and 
$\epsilon=0.005$ (black).  
The view windows are zoomed into the regions where different $\epsilon$ values impose 
different solution profiles. The results are discussed in example 11. }
\label{diff_epsilon} \end{figure}

\noindent{\bf Example 11} $\epsilon = 0.001, \epsilon = 0.002, \epsilon = 0.003, \epsilon = 0.004, \epsilon = 
0.005$.\\
In this example, we will compare the solution profiles for different $\epsilon$ values. Fixing $T=0.5, \Delta x = 
0.0001, \lambda = \frac{\Delta t}{\Delta x} = 0.1$, we show the numerical results in Figure \ref{diff_epsilon} 
for $\epsilon= 0.001$ (blue), $\epsilon=0.002$ (yellow), $\epsilon=0.003$ (magenta),  $\epsilon=0.004$ (green), 
and $\epsilon=0.005$ (black). 
For the purpose of cross reference, we choose the same nine sets of parameter settings as in examples 1-- 6. To 
assist the observation, the figures in Figure \ref{diff_epsilon}  are zoomed into the regions where different 
$\epsilon$ values introduce different solution profiles.
The numerical solutions clearly show that as $\epsilon$ increases, the numerical solution is smeared out, and the 
jump location becomes less accurate. 
Notice that $\tau$ is responsible for the competition between the diffusion and dispersion, which in turn 
determines the plateau values. Hence varying $\epsilon$ value doesn't affect the plateau location.
\\

\section{Conclusion}
\label{conclusion}
We proved that the solution to the infinite domain problem can be approximated by that of the bounded domain problem. This provides a theoretical justification for using finite domain to calculation the numerical solution of the MBL equation (\ref{MBL}). We also extended the classical central scheme originally designed for the hyperbolic systems to solve the MBL equation, which is of pseudo-parabolic type.
The numerical solutions for qualitatively different parameter values $\tau$ and initial conditions $u_B$ show that the jump locations are consistent with the theoretical calculation and the plateau heights are consistent with the numerically obtained values given in \cite{duijn}.
In particular, when $\tau>\tau_*$, for $u_B\in(\underline{u},\bar{u})$, the numerical solutions give  non-monotone water saturation profiles, which is consistent with the experimental observations.
In addition, the order tests show that the proposed second and third order central schemes achieved the desired accuracies. 

In \cite{Ying,duijn2}, the two-dimensional space extension of the modified Buckley-Leverett equation has been derived. One of the 
future directions is to develop high order numerical schemes to solve the two-dimensional MBL equation.
Central schemes have been used to solve high dimensional hyperbolic problem and dispersive problem (\cite{Jiang,Levy}), which makes
it a good candidate for such a task.

\begin{appendix}

\section{Proof of the lemmas}
\label{proof_lemma}

\begin{proof}[Proof to lemma \ref{lemmafu}]
Let $g(u)=\frac{f(u)}{u}=\frac{u}{u^2+M(1-u)^2}$, then
\begin{eqnarray*}
g'(u)=\frac{M-(1+M)u^2}{(u^2+M(1-u)^2)^2}\left\{
\begin{array}{lll}
>0&\mathrm{if}&0<u<\sqrt{\frac{M}{M+1}}\\
=0&\mathrm{if}&u=\sqrt{\frac{M}{M+1}}\\
<0&\mathrm{if}&u>\sqrt{\frac{M}{M+1}}
\end{array}
\right.
\end{eqnarray*}
and hence $g(u)$ achieves its maximum at $u=\sqrt{\frac{M}{M+1}}$. Therefore, $\frac{f(u)}{u}=g(u)\le D$, where $D=\frac{f(\alpha)}{\alpha}$ and $\alpha=\sqrt{\frac{M}{M+1}}$, and in turn, we have that $f(u)\le Du$ for all $0\le u \le 1$.
\end{proof}

\begin{proof}[Proof to lemma \ref{lemma_old1} (\ref{lemma1})]
\begin{eqnarray*}
\int_{0}^{+\infty}
\left|e^{-\frac{x+\xi}{\epsilon\sqrt{\tau}}}-e^{-\frac{|x-\xi|}{\epsilon\sqrt{\tau}}}\right|
e^{\frac{\lambda x-\lambda\xi}{\epsilon\sqrt{\tau}}}\,d\xi
%
%
%
%
%
%
%
=\epsilon\sqrt{\tau}\frac{-2+2e^{\frac{(\lambda-1)x}{\epsilon\sqrt{\tau}}}}{\lambda^2-1}
\le\frac{2\epsilon\sqrt{\tau}}{1-\lambda^2} \quad\text{  if }\lambda\in(0,1).
\end{eqnarray*}
\end{proof}

\begin{proof}[Proof to lemma \ref{lemma_old1} (\ref{lemma2})]
\begin{eqnarray*}
\int_{0}^{+\infty}
\left|e^{-\frac{x+\xi}{\epsilon\sqrt{\tau}}}-e^{-\frac{|x-\xi|}{\epsilon\sqrt{\tau}}}\right|
e^{\frac{\lambda x-\xi}{\epsilon\sqrt{\tau}}}\,d\xi
%
%
%
%
%
%
=xe^{\frac{(\lambda-1)x}{\epsilon\sqrt{\tau}}} \le
\frac{\epsilon\sqrt{\tau}}{e(1-\lambda)} \quad\text{  if } \lambda\in(0,1).
\end{eqnarray*}
\end{proof}

\begin{proof}[Proof to lemma \ref{lemma_old1} (\ref{lemma3})]
Based on the assumption on $u_0$ in (\ref{u0})
\begin{equation}
\begin{split}
\int_{0}^{+\infty}
\left|e^{-\frac{x+\xi}{\epsilon\sqrt{\tau}}}-e^{-\frac{|x-\xi|}{\epsilon\sqrt{\tau}}}\right|e^{\frac{\lambda
x}{\epsilon\sqrt{\tau}}} |u_0(\xi)|\,d\xi
\le& \int_{0}^{+\infty}
e^{-\frac{|x-\xi|}{\epsilon\sqrt{\tau}}}e^{\frac{\lambda
x}{\epsilon\sqrt{\tau}}} |u_0(\xi)|\,d\xi  \\
\le&
C_ue^{\frac{\lambda
x}{\epsilon\sqrt{\tau}}}\int_{0}^{L_0}e^{-\frac{|x-\xi|}{\epsilon\sqrt{\tau}}}
\,d\xi
%
%
=C_u y_1(x) 
\end{split}
\label{u01}
\end{equation}
Calculating $y_1(x)$ with the assumption that $\lambda\in(0,1)$,
we get
\begin{eqnarray*}
y_1(x)=\left\{
\begin{array}{lll}
e^{\frac{\lambda
x}{\epsilon\sqrt{\tau}}}\int_0^{L_0}e^{-\frac{|x-\xi|}{\epsilon\sqrt{\tau}}}\,d\xi\le2\epsilon\sqrt{\tau}e^{\frac{\lambda
x}{\epsilon\sqrt{\tau}}}\le 2\epsilon\sqrt{\tau}e^{\frac{\lambda
L_0}{\epsilon\sqrt{\tau}}} & & \mathrm{for}~~x\in[0,L_0]\\
&&\\
e^{\frac{(\lambda-1)x}{\epsilon\sqrt{\tau}}}\int_0^{L_0}e^{\frac{\xi}{\epsilon\sqrt{\tau}}}\,d\xi\le\epsilon\sqrt{\tau}e^{\frac{(\lambda-1)x+L_0}{\epsilon\sqrt{\tau}}}
\le\epsilon\sqrt{\tau}e^{\frac{\lambda L_0}{\epsilon\sqrt{\tau}}}
&&\mathrm{for}~~x\in[L_0,+\infty)
\end{array} \right.
\end{eqnarray*}
Therefore, we get the desired inequality
\begin{eqnarray*}
\int_{0}^{+\infty}\left|e^{-\frac{x+\xi}{\epsilon\sqrt{\tau}}}-e^{-\frac{|x-\xi|}{\epsilon\sqrt{\tau}}}\right|e^{\frac{\lambda
x}{\epsilon\sqrt{\tau}}} |u_0(\xi)|\,d\xi \le
2C_u\epsilon\sqrt{\tau}e^{\frac{\lambda
L_0}{\epsilon\sqrt{\tau}}}
.
\end{eqnarray*}
\end{proof}

\begin{proof}[Proof to lemma  \ref{lemma_old2} (\ref{lemma11})]$ $\\
\begin{eqnarray*}
&&\int_{0}^{+\infty}
\left|e^{-\frac{x+\xi}{\epsilon\sqrt{\tau}}}+\mathrm{sgn}(x-\xi)e^{-\frac{|x-\xi|}{\epsilon\sqrt{\tau}}}\right|
e^{\frac{\lambda x-\lambda\xi}{\epsilon\sqrt{\tau}}}\,d\xi
\\
&=&
\frac{\epsilon\sqrt{\tau}}{\lambda^2-1}
\left(-2+2\lambda
e^{\frac{(\lambda-1)x}{\epsilon\sqrt{\tau}}}-2(\lambda-1)e^{-\frac{2x}{\epsilon\sqrt{\tau}}}\right)
\le\frac{2\epsilon\sqrt{\tau}}{1-\lambda^2}\qquad\text{ if } \lambda\in(0,1).
\end{eqnarray*}
\end{proof}

\begin{proof}[Proof to lemma  \ref{lemma_old2} (\ref{lemma12})]$ $\\
\begin{eqnarray*}
&&\int_{0}^{+\infty}
\left|e^{-\frac{x+\xi}{\epsilon\sqrt{\tau}}}+\mathrm{sgn}(x-\xi)e^{-\frac{|x-\xi|}{\epsilon\sqrt{\tau}}}\right|
e^{\frac{\lambda x-\xi}{\epsilon\sqrt{\tau}}}\,d\xi\\
%
%
%
%
%
%
%
%
%
%
&=&
\frac{2e^{{\frac{(\lambda-3)x}{\epsilon\sqrt{\tau}}}}-2e^{{\frac{(\lambda-1)x}{\epsilon\sqrt{\tau}}}}}{\frac{-2}{\epsilon\sqrt{\tau}}}
+xe^{\frac{(\lambda-1)x}{\epsilon\sqrt{\tau}}}
\le \epsilon\sqrt{\tau}+\frac{\epsilon\sqrt{\tau}}{e(1- \lambda)}\qquad\text{ if }\lambda\in(0,1).
\end{eqnarray*}
\end{proof}

\begin{proof}[Proof to lemma \ref{lemma_old2} (\ref{lemma13})]
Based on the assumption on $u_0$ in (\ref{u0})
\begin{equation}
\begin{split}
&\int_{0}^{+\infty}
\left|e^{-\frac{x+\xi}{\epsilon\sqrt{\tau}}}+\mathrm{sgn}(x-\xi)e^{-\frac{|x-\xi|}{\epsilon\sqrt{\tau}}}\right|
e^{\frac{\lambda x}{\epsilon\sqrt{\tau}}}|u_0(\xi)|\,d\xi  \\
\le&
~C_ue^{\frac{\lambda
x}{\epsilon\sqrt{\tau}}}\int_{0}^{L_0}\left|e^{-\frac{x+\xi}{\epsilon\sqrt{\tau}}}+\mathrm{sgn}(x-\xi)
e^{-\frac{|x-\xi|}{\epsilon\sqrt{\tau}}}\right|\,d\xi
\\
=& ~C_u y_3(x)
\end{split}
\label{u02}
\end{equation}
Calculating $y_3(x)$ with the assumption that $\lambda\in(0,1)$,
we get for $x\in[0,L_0]$
\begin{eqnarray*}
y_3(x)\le
e^{\frac{(\lambda-1)x}{\epsilon\sqrt{\tau}}}\int_0^x(e^{-\frac{\xi}{\epsilon\sqrt{\tau}}}+e^{\frac{\xi}{\epsilon\sqrt{\tau}}})\,d\xi+e^{\frac{(\lambda+1)x}{\epsilon\sqrt{\tau}}}\int_x^{L_0}e^{-\frac{\xi}{\epsilon\sqrt{\tau}}}\,d\xi\le2\epsilon\sqrt{\tau}e^{\frac{\lambda L_0}{\epsilon\sqrt{\tau}}}
\end{eqnarray*}
and 
\begin{eqnarray*}
y_3(x)\le
e^{\frac{(\lambda-1)x}{\epsilon\sqrt{\tau}}}\int_0^{L_0}(e^{-\frac{\xi}{\epsilon\sqrt{\tau}}}+e^{\frac{\xi}{\epsilon\sqrt{\tau}}})\,d\xi\le\epsilon\sqrt{\tau}e^{\frac{(\lambda-1)x+L_0}{\epsilon\sqrt{\tau}}}
\le\epsilon\sqrt{\tau}e^{\frac{\lambda L_0}{\epsilon\sqrt{\tau}}}
\end{eqnarray*}
for $x\in[L_0,+\infty)$.

Therefore, we get the desired inequality
\begin{eqnarray*}
\int_{0}^{+\infty}
\left|e^{-\frac{x+\xi}{\epsilon\sqrt{\tau}}}+\mathrm{sgn}(x-\xi)e^{-\frac{|x-\xi|}{\epsilon\sqrt{\tau}}}\right|
e^{\frac{\lambda x}{\epsilon\sqrt{\tau}}}|u_0(\xi)|\,d\xi\le
2C_u\epsilon\sqrt{\tau}e^{\frac{\lambda L_0}{\epsilon\sqrt{\tau}}}
.
\end{eqnarray*}

\end{proof}

\begin{proof}[Proof to lemma \ref{lemma_old4} (\ref{lemmaphi1})]
\begin{eqnarray*}
\left|\phi_1(x)-e^{-\frac{x}{\epsilon\sqrt{\tau}}}\right|
=e^{-\frac{L}{\epsilon\sqrt{\tau}}}\left|\frac{e^{-\frac{x}{\epsilon\sqrt{\tau}}}-e^{\frac{x}{\epsilon\sqrt{\tau}}}}{e^{\frac{L}{\epsilon\sqrt{\tau}}}-e^{-\frac{L}{\epsilon\sqrt{\tau}}}}\right|
=e^{-\frac{L}{\epsilon\sqrt{\tau}}}\left|\phi_2(x)\right| .
\end{eqnarray*}
\end{proof}

\begin{proof}[Proof to lemma \ref{lemma_old4} (\ref{lemmaphi2})]
Since
$
\phi_2(x)=
\frac{e^{\frac{x}{\epsilon\sqrt{\tau}}}-e^{-\frac{x}{\epsilon\sqrt{\tau}}}}
{e^{\frac{L}{\epsilon\sqrt{\tau}}}-e^{-\frac{L}{\epsilon\sqrt{\tau}}}}
$, we see that
$
\phi_2'(x)=\frac{1}{\epsilon\sqrt{\tau}}\frac{e^{\frac{x}{\epsilon\sqrt{\tau}}}+e^{-\frac{x}{\epsilon\sqrt{\tau}}}}{e^{\frac{L}{\epsilon\sqrt{\tau}}}-e^{-\frac{L}{\epsilon\sqrt{\tau}}}}>0
$
and hence
$\phi_2(x)\le\phi_2(L)=1$  for $x\in[0,L]$.
\end{proof}

\begin{proof}[Proof to lemma  \ref{lemma_old4} (\ref{lemmaphi2prime})]
$\phi_2'(x)=\frac{1}{\epsilon\sqrt{\tau}}\frac{e^{\frac{x}{\epsilon\sqrt{\tau}}}+e^{-\frac{x}{\epsilon\sqrt{\tau}}}}{e^{\frac{L}{\epsilon\sqrt{\tau}}}-e^{-\frac{L}{\epsilon\sqrt{\tau}}}}$ gives that $\phi_2''(x)=\frac{1}{\epsilon^2\tau}\phi_2(x)>0$, and hence $\phi_2'(x)\le\phi_2'(L)=\frac{1}{\epsilon\sqrt{\tau}}\frac{e^{\frac{L}{\epsilon\sqrt{\tau}}}+e^{-\frac{L}{\epsilon\sqrt{\tau}}}}{e^{\frac{L}{\epsilon\sqrt{\tau}}}-e^{-\frac{L}{\epsilon\sqrt{\tau}}}}=\frac{1}{\epsilon\sqrt{\tau}}\frac{e^{\frac{2L}{\epsilon\sqrt{\tau}}}+1}{e^{\frac{2L}{\epsilon\sqrt{\tau}}}-1}\le\frac{2}{\epsilon\sqrt{\tau}}$ if $\epsilon\ll 1$ for $x\in[0,L]$.
\end{proof}


\end{appendix}

\section*{Acknowledgments}
CYK would like to thank Prof. L.A. Peletier for introducing MBL equation
and Mathematical Biosciences Institute at OSU for the hospitality and
support.
%
%

\bibliography{1D}
\bibliographystyle{amsplain}

\end{document}